\documentclass[a4paper,DIV=12,10pt]{scrartcl}

\usepackage[utf8]{inputenc}
\usepackage{amsfonts,amssymb,amsmath}
\usepackage{epsfig}
\usepackage{MnSymbol}
\usepackage{booktabs}
\usepackage{cite}
\usepackage{float}
\usepackage{hyperref}

\usepackage{graphics}
\usepackage{subfig}
\usepackage{color}

\usepackage{enumitem}

\newcommand{\ud}{\,\mathrm{d}}
\newcommand{\R}{\mathbb{R}}
\newcommand{\N}{\mathbb{N}}

\renewcommand{\vec}{\boldsymbol}

\newcommand{\grd}{\operatorname{grad}_0}
\newcommand{\Grd}{\operatorname{Grad}_0}
\newcommand{\dv}{\operatorname{div}}
\newcommand{\Dv}{\operatorname{Div}}
\newcommand{\dggrd}{\operatorname{grad}_{\operatorname{dg}}}
\newcommand{\dgdiv}{\operatorname{div}_{\operatorname{dg}}}
\newcommand{\dgGrd}{\operatorname{Grad}_{\operatorname{dg}}}
\newcommand{\dgDiv}{\operatorname{Div}_{\operatorname{dg}}}

\newcommand{\Hhnc}{H_h}
\newcommand{\hItau}{\widehat I_\tau^{\,k+1}}
\newcommand{\InGR}{I_\tau^{n}}

\newcommand{\olq}{\overline{q}}

\newcommand*{\ldblbrace}{\{\mskip-5mu\{}
\newcommand*{\rdblbrace}{\}\mskip-5mu\}}

\numberwithin{equation}{section}
\numberwithin{figure}{section}
\numberwithin{table}{section}

\newtheorem{defi}{Definition}[section]
\newtheorem{thm}[defi]{Theorem}
\newtheorem{lem}[defi]{Lemma}
\newtheorem{rem}[defi]{Remark}
\newtheorem{cor}[defi]{Corollary}

\newtheorem{prob}[defi]{Problem}

\newenvironment{mproof}{\paragraph{Proof.}}{\hfill$\blacksquare$ \medskip}

\begin{document}
	
\title{Structure preserving discontinuous Galerkin approximation of a hyperbolic-parabolic system}
	
\author{
Markus Bause$^\dag$\thanks{bause@hsu-hh.de (corresponding author)}\;, Sebastian Franz$^\ddag$\\
{\small ${}^\dag$ Helmut Schmidt University, Faculty of Mechanical and Civil Engineering, Holstenhofweg 85,}\\ 
{\small 22043 Hamburg, Germany}\\
{\small $^\ddag$ Technical University of Dresden, Institute of Scientific Computing, 01062 Dresden, Germany}\\
{\small Germany}
}
	
	
	\maketitle
	
\begin{abstract}
\textbf{Abstract.}
We study the numerical approximation of a coupled hyperbolic-parabolic system by a family of discontinuous Galerkin space-time finite element methods. The model is rewritten as a first-order evolutionary problem that is treated by the unified abstract solution theory of R.\ Picard. For the discretization in space, generalizations of the distribution gradient and divergence operators on broken polynomial spaces are defined. Since their skew-selfadjointness is perturbed by boundary surface integrals, adjustments are introduced such that the skew-selfadjointness of the first-order differential operator in space is recovered. Well-posedness of the fully discrete problem and error estimates for the discontinuous Galerkin approximation in space and time are proved. 
\end{abstract}
	
	
	
\section{Introduction}
\label{Sec:Introduction}

We study the numerical approximation by discontinuous Galerkin methods in space and time of solutions to the hyperbolic-parabolic system 
\begin{subequations}
	\label{Eq:HPS}
	\begin{alignat}{3}
		\label{Eq:HPS_1}
		\rho \partial_t^2 u - \nabla \cdot (C \varepsilon (u)) + \alpha \nabla p & = \rho  f\,, && \quad \text{in } \;
		\Omega  \times (0,T]\,,\\[1ex]
		\label{Eq:HPS_2}
		c_0\partial_t p + \alpha  \nabla \cdot \partial_t  u  -  \nabla \cdot (K \nabla p)  & = g\,, && \quad \text{in } \; \Omega \times (0,T]\,,\\[1ex]
		\label{Eq:HPS_3}
		u (0) = u_0\,, \quad \partial_t  u (0) = u_1\,, \quad p(0) & = p_0\,, && \quad \text{in } \; \Omega \,,\\[1ex]
		\label{Eq:HPS_4}
		u = 0\,, \quad p& =0\,,  && \quad \text{on } \; \partial\Omega \times (0,T]\,. 
	\end{alignat}
\end{subequations}
For this, we rewrite \eqref{Eq:HPS} as a first order evolutionary problem in space and time on the open bounded domain $\Omega \subset \R^d$, with $d\in \{2,3\}$, and for the final time $T>0$. System \eqref{Eq:HPS} is investigated as a prototype model problem for poro- and thermoelasticity; cf., e.g., \cite{B41,B55,B72,JR18,S00}. In poroelasticity, Eqs.\ \eqref{Eq:HPS_1} and \eqref{Eq:HPS_2} describe the conservation of momentum and mass. The unknowns are the effective solid phase displacement $u$ and the effective fluid pressure $p$. The quantity $\varepsilon (u):= ( \nabla u + (\nabla u)^\top)/2$ denotes the symmetrized gradient or strain tensor. Further, $\rho$ is the effective mass density, $C$ is Gassmann’s fourth order effective elasticity tensor, $\alpha$ is Biot’s pressure-storage coupling tensor, $c_0$ is the specific storage coefficient and $K$ is the permeability field. For simplicity, the positive quantities $\rho>0$, $\alpha>0$ and $c_0 >0$ are assumed to be constant in space and time. Moreover, the tensors $C$ and $K$ are assumed to be symmetric and positive definite and independent of the space and time variables as well. In \eqref{Eq:HPS_1}, the effects of secondary consolidation (cf.~\cite{MC96}), described in certain models by the additional term $\lambda^\ast \delta_{ij} \varepsilon (\partial_t u)$ in the total stress, are not included here. Beyond the classical applications of \eqref{Eq:HPS} in subsurface hydrology and geophysics, for instance in reservoir engineering, systems like \eqref{Eq:HPS} have recently attracted reseachers' interest in biomedical engineering; cf., e.g., \cite{C16,CADQ23,N13}. In thermoelasticity, the system \eqref{Eq:HPS} describes the flow of heat through an elastic structure. In that context, $p$ denotes the temperature, $c_0$ is the specific heat of the medium, and $K$ is the conductivity. The homogeneous Dirichlet boundary conditions \eqref{Eq:HPS_4} are studied here for simplicity and brevity.

By introducing the variable of unknowns $U = (v, \sigma,p,\olq)^\top$, with the quantities $v:=\partial_t u$, $\sigma:=C \varepsilon$ and $\overline q:=-K\nabla p+\alpha v$, we transform the system \eqref{Eq:HPS} into an abstract evolutionary equation written as the sum of two unbounded first-order differential operators, one of them involving a first order differential operator in time and the other one involving first order differential operators in space. In the exponentially weighted in time Bochner space $H_\nu(\R;H)$ defined in \eqref{Eq:Hnu}, with some weight $\nu>0$ and the Hilbert space $H=L^2(\Omega)^{(d+1)^2}$, we then obtain the evolutionary equation for $U$ that 
\begin{equation}
\label{Eq:EP}
(\partial_t M_0 + M_1 + A)  U = F\,. 
\end{equation}
In \eqref{Eq:EP}, $M_0$ and $M_1$ are bounded linear selfadjoint operators in $H$ and $A$ is an unbounded skew-selfadjoint operator in $H$. The right-hand side function $F$ in \eqref{Eq:EP} depends on the source terms $f$ and $g$ of \eqref{Eq:HPS}. For \eqref{Eq:EP}, a solution mechanism developed by R.\ Picard \cite{P09} can be applied. It is based on monotonicity of both the sum of the mentioned unbounded operators together with its adjoint computed in the space-time Hilbert space. For the presentation of the solution theory we refer to \cite[Thm.~6.2.1]{STW22}. The well-posedness criterion for \eqref{Eq:EP}, that is summarized in Thm.~\ref{Thm:WP}, is elementary and general. It can be verified without dealing with the intricacies of more involved solution methods. This is an appreciable advantage of Picard's theorem \cite{P09}. A priori, there is no explicit initial condition implemented in the theory. For \eqref{Eq:EP}, an initial condition of the form $\lim_{t\searrow t_0}M_0 U(t)=M_0 U_0$ for some $t_0\in \R$ and $U\in D(A)$ can be implemented by a distributional right-hand side term $F + \delta_{t_0} M_0 U_0$ for some $F\in H_\nu(\R;H)$ supported on $[t_0,\infty)$ and the Dirac distribution $ \delta_{t_0}$ at $t_0$. For details of this we refer to \cite[Sec.~6.2.5]{PMcG11} and \cite[Chap.~9]{STW22}. By introducing the four-field formulation for the unknown vector $U$ the problem size is increased. However, in poroelasticity the explicit approximation of the flux variable $q=-K\nabla p$ is often desirable and of higher importance than the approximation of the fluid pressure itself. For instance, this holds if reactive transport of species, dissolved in the fluid, is studied further. Simulations then demand for accurate approximations of the flux variable $q$. A similar argument applies to the stress tensor $\sigma$, if this variable is the goal quantity of physical interest in \eqref{Eq:HPS} or needs to be post-processed for elucidating phenomena modeled by \eqref{Eq:HPS}. In implementations, the symmetrie of the stress tensor $\sigma$ can still be exploited to reduce the problem's size.  
 
In this work we propose and analyze fully discrete numerical approximation schemes that are built for the evolutionary equation \eqref{Eq:EP}. Their key feature is that they essentially preserve the abstract evolutionary form \eqref{Eq:EP} and the operators' properties. However, due to a nonconforming discretization in space that is applied here, the skew-selfadjointness of the discrete counterpart of $A$ in \eqref{Eq:EP} is perturbed by non-vanishing contributions arising from boundary face integrals. Therefore,  a correction term is introduced on the discrete level to overcome this defect and ensure that a discrete counterpart of the skew-selfadointness that is essentially used in the analyses is satisfied. In the design of numerical methods, structure preserving approaches ensuring that important properties of differential operators and solutions to the continuous problem are maintained on the fully discrete level are highly desirable and important to ensure physical realism of numerical predictions. We focus discontinuous Galerkin (DG) discretizations of the space and time variables. DG methods for the space discretization (cf., e.g., \cite{PE12,DF15,R08}) have shown their high flexibility and accuracy in approximating reliably solutions to partial differental equations, even solutions with complex structures or discontinuities and in anisotropic or heterogenous media. The application of DG schemes for the space discretization of  \eqref{Eq:EP} and the definition of the DG counterpart of $A$ in  \eqref{Eq:EP}  to preserve skew-selfadjointness represent the key innovation of this work over a series of previous ones \cite{ABMS23,FTW19,F21,F23} based on Picard's theory. For the DG space discretization, the definition of the distribution gradient and divergence operator is extended to broken polynomial spaces by penalizing the jumps of the unknowns over interelement surfaces. By still adding some boundary correction due to the nonconformity of DG methods, the skew-selfadjointness of $A$ is passed on to its discrete counterpart $A_h$.  This consistent definition and treatment of the DG gradient and DG divergence operators for the nonconforming approximation is essential for the overall approach and its analysis. It has not been studied yet. 

For the discretization in time we use the DG method \cite{T06}. Variational time discretizations offer the appreciable advantage of the natural construction of families of schemes with higher order members, even for complex coupled systems of equations. There exists a strong link to Runge--Kutta methods; cf.\ \cite{AMN09,AMN11}. DG time discretizations are known to be strongly A-stable. For elastodynamics and wave propagation they violate the energy conservation principle of solutions to the continuous problem. This might evoke effects of damping or dispersion. However, the convergence of the jump terms at the discrete time nodes  is ensured; cf.\ \cite[Thm.\ 2.3]{F23}. Continuous in time Galerkin (CG) methods (cf.\, e.g., \cite{ABBM20,BKR22,BKRS20,F21,AM04} and the references therein) are known to be A-stable only, but they preserve the energy of solutions \cite[Sec.~6]{BKRS20}. These families are more difficult to analyse since they lead to Galerkin--Petrov methods with trial and test spaces differing from each other. For this reason and due to computational advantages gained for simulations of the second-order form \eqref{Eq:HPS}, DG time discretizations are studied here. For studies of CG schemes with continuous in time discrete solutions we refer to, e.g., \cite{ABBM20,BKR22,BKRS20,F21,KM04,DFW16} and the references therein. For a numerical study of DG and CG time discretizations of \eqref{Eq:HPS} we refer to \cite{ABMS23_2}.

In \cite{FTW19} and \cite{F21}, one of the authors of this work studies with his coauthors numerical schemes based on DG and CG Galerkin methods in time and conforming Galerkin methods in space for evolutionary problems \eqref{Eq:EP} of changing type. By decomposing $\Omega$ into three disjoint sets and defining the $M_0$ and $M_1$ setwise, the system \eqref{Eq:EP} degenerates to elliptic, parabolic or hyperbolic type on these sets. Usually, degenerating problems are difficult to analyze. Due to the weak assumptions about the operators made in the theory of Picard \cite{P09}, such type of problems can be embedded into this framework. The same applies to the concept of perfectly matched layers in wave propagation; cf., e.g. \cite{B94,CM98}. They are used to truncate the entire space $\R^d$ or unbounded domains to bounded computational ones and mimic non-reflecting boundary conditions. The analysis of wave propagation with artifical absorbing layer and changing equations in either regions becomes feasible as well by the abstract solution theory. 

In \cite{CDW22}, space-time DG methods for weak solutions of hyperbolic linear first-order symmetric Friedrichs systems describing acoustic, elastic, or electro-magnetic waves are proposed. For an introduction into the theory of first-order symmetric Friedrichs systems we refer to \cite{EG06a,EG06b,EG08} and \cite[Chap.~7]{PE12}. Similarly to this work, in \cite{CDW22} a first-order in space and time formulation of a second-order hyperbolic problem is used. In contrast to this work, no coupled system of mixed hyperbolic-parabolic type is considered there. In \cite{CDW22}, the mathematical tools for proving well-posedness of the space-time DG discretization and error estimates are based on the theory of first-order Friedrichs systems. The theory strongly differs from Picard's theorem \cite{P09} that is used here. The differences of either approaches still require elucidation. In deriving space-time DG methods and proving error estimates, differences become apparent in the norms with respect to that convergence is proved. In \cite{CDW22}, stability and convergence estimates are provided with respect to a mesh-dependent DG norm that includes the $L^2$ norm at the final time; cf.\ also \cite{BMPS21}. 

Here, convergence of the fully discrete approximation $U_{\tau,h}$ of \eqref{Eq:EP} is proved in Thm.~\ref{Thm:ErrEstSTnc} with respect to the natural and induced norm of the exponentially weighted Bochner space $H_\nu(\R;H)$, with the product space $H=L^2(\Omega)^{(d+1)^2}$ equipped with the $L^2$-norm. For the full discretization $U_{\tau,h}$ of the solution $U$ to \eqref{Eq:EP} we show that
\begin{equation}
	\label{KeyRes}
 	\sup_{t\in [0,T]} \langle M_0 (U-U_{\tau,h}))(t), (U-U_{\tau,h}))(t) \rangle_H +  \| U- U_{\tau,h}\|^2_{\nu} \leq C (\tau^{2(k+1)} + h^{2r})\,,
\end{equation}
where $\|\cdot \|_{\nu}$ is the exponentially weighted natural norm associated with $H_\nu(\R;H)$. Further, $k$ and $r$ are the piecewise polynomial degrees in time and space, respectively. 

The paper is organized as follows. In Sec.~\ref{Sec:EvFoWP} the evolutionary form \eqref{Eq:EP} of \eqref{Eq:HPS} is derived and its well-posedness is shown. The space-time discretization of \eqref{Eq:EP} by the DG method is presented in Sec.~\ref{Sec:TimeDisc}. Its error analysis is addressed in Sec.~\ref{Subsec:ErrEst}. In Sec.~\ref{Sec:SumOut}, we end with a summary and outlook. 

\section{Evolutionary formulation and its well-posedness}
\label{Sec:EvFoWP}

In this section we rewrite formally the coupled hyperblic-parabolic problem \eqref{Eq:HPS} as an evolutionary problem \eqref{Eq:EP} by introducing auxiliary variables. For the  evolutionary problem we present a result of well-posedness that is based on Picard's theorem; cf.\ \cite{P09} and \cite[Thm.~6.2.1]{STW22}. Therein, the evolutionary problem is investigated on the whole time axis, for  $t\in \R$, in the exponentially weighted Bochner space $H_\nu(\R;H)$ introduced in Def.~\ref{Def:Hnu}. Throughout, we use usual notation for standard Sobolev spaces. In the notation, we do not differ between scalar-, vector- or tensor-valued functions. 

\begin{defi}
\label{Def:Hnu}
Let $H$ be a real Hilbert space with associated norm $\|\cdot \|_H$. For $\nu >0$, we put
\begin{equation}
\label{Eq:Hnu}
H_\nu(\R;H) := \Big\{f:\R \rightarrow H  : \int_{\R} \| f(t) \|_H^2 \operatorname{e}^{-2\nu t} \ud t \Big\}\,.
\end{equation}		
\end{defi}	 

The space $H_\nu(\R;H)$, equipped with the inner product 
\begin{equation}
\label{Def:SPHnu}
\langle f,g\rangle_{\nu} := \int_{\R} \langle f(t),g(t)\rangle_H \operatorname{e}^{-2\nu t} \ud t\,, \quad \text{for}\;  f,g\in H_\nu(\R;H)\,, 
\end{equation}
is a Hilbert space. The norm induced by the inner product \eqref{Def:SPHnu} is denoted by $\|\cdot \|_\nu$. Moreover, we define $\partial_t$ to be the closure of the operator 
\begin{equation*}
	\partial_t : C^\infty_c (\R;H) \subset  H_\nu(\R;H)  \rightarrow H_\nu(\R;H)\,, \;\; \phi \rightarrow \phi^\prime \,,
 \end{equation*}
where $C^\infty_c (\R;H) $ is the space of infinitely differentiable $H$-valued functions on $\R$ with compact support. The domain of the time derivative of $n$-order, denoted by $\partial_t^s$, is the space $H^s_\nu(\R;H)$. Before rewriting  \eqref{Eq:HPS} in the form \eqref{Eq:EP}, we need to introduce differential operators with respect to the spatial variables. 

\begin{defi}
Let $\Omega \subset \R^d$, for $d\in \N$, be an open non-empty set.	 Then we define 
\begin{equation*}
L^2(\Omega)^{d\times d}_{\operatorname{sym}} := \left\{(\phi_{ij})_{i,j=1,\ldots,d}\in L^2(\Omega) : \phi_{ij} = \phi_{ji} \;\; \forall \; i,j\in \{1,\ldots,d\}\right\}	\,.
\end{equation*}
\end{defi}

\begin{defi}
Let $\Omega \subset \R^d$, for $d\in \N$, be an open non-empty set.	 We put 
\begin{equation}
\label{Def:grad0}
\grd : H^1_0(\Omega) \subset L^2(\Omega) \rightarrow L^2(\Omega)^d\,,  \;\;  \phi \rightarrow (\partial_j \phi)_{j=1,\ldots,d}\,,
\end{equation}
and 
\begin{equation}
\label{Def:Grad0} 
\Grd : H^1_0(\Omega)^d \subset L^2(\Omega)^d \rightarrow L^2(\Omega)^{d\times d}_{\operatorname{sym}}\,,   \;\;  (\phi_j)_{j=1,\ldots,d}  \rightarrow \dfrac{1}{2} (\partial_l\phi_j + \partial_j \phi_l)_{j,l=1,\ldots ,d}\,.
\end{equation}
Moreover, we put 
\begin{equation}
	\label{Def:div}
	\dv : D(\dv) \subset L^2(\Omega)^d \rightarrow L^2(\Omega)\,,  \;\;   \dv := - (\grd)^\ast\,, 
\end{equation}
and 
\begin{equation}
	\label{Def:Div}
	\Dv : D(\Dv) \subset L^2(\Omega)^{d\times d}_{\operatorname{sym}} \rightarrow L^2(\Omega)^d\,,  \;\;   \Dv := - (\Grd)^\ast\,.
\end{equation}
\end{defi}	

We note that $\Grd u = \varepsilon(u)$ for $u\in H^1_0(\Omega)^d$. The operator $\dv$ in \eqref{Def:div} assigns each $L^2$ vector field its distributional divergence with maximal domain, that is, 
\begin{equation*}
D(\dv) = 	\left\{\phi \in L^2(\Omega)^d : \sum_{i=1}^d \partial_i \phi_i\in L^2(\Omega)\right\}\,.
\end{equation*}
Similarly, the operator $\Dv$ in \eqref{Def:Div} assigns each  $L^2(\Omega)^{d\times d}_{\operatorname{sym}}$  tensor field its distributional divergence with maximal domain, that is, 
\begin{equation*}
	D(\Dv) = 	\left\{\phi \in L^2(\Omega)^{d\times d}_{\operatorname{sym}}: \left(\sum_{i=1}^d \partial_i \phi_{ij}\right)_{j=1,\ldots,d} \in L^2(\Omega)^d\right\}\,.
\end{equation*}

To rewrite \eqref{Eq:HPS} formally as a first-order evolutionary problem, we introduce the set of new unknowns
\begin{equation}
\label{Def:NU}
v:= \partial_t u\,, \quad \sigma := C\varepsilon \quad \text{and} \quad   q:= -K \nabla p\,.
\end{equation}
Using \eqref{Def:NU} and differentiating the second of the definitions in \eqref{Def:NU} with respect to the time variable, we recast \eqref{Eq:HPS_1} and  \eqref{Eq:HPS_2} as the first order in space and time system  
\begin{subequations}
	\label{Eq:HPS_10}
	\begin{alignat}{3}
		\label{Eq:HPS_11}
		\rho \partial_t v - \Dv \sigma + \alpha \grd p & = \rho  f\,,\\[1ex]
		\label{Eq:HPS_12}
		S \partial_t \sigma - \Grd v & = 0 \,, \\[1ex]
		\label{Eq:HPS_13}
		c_0\partial_t p + \alpha \dv  v + \dv q & = g\,, \\[1ex]
		\label{Eq:HPS_14}
	    K^{-1} q +  \grd p & = 0 \,,
	\end{alignat}
\end{subequations}
where $S$ denotes the positive definite, fourth order compliance tensor of the inverse stress-strain relation of Hook's law of linear elasticity, 
\begin{equation}
 \varepsilon = S \sigma\,.
\end{equation}
In matrix-vector notation the system \eqref{Eq:HPS_10} reads as  
\begin{equation}
\label{Eq:HPS_20}
\left(\partial_t 
\begin{pmatrix}
	\rho & 0 & 0 &0\\
	0 & S & 0 & 0 \\
	0 & 0 & c_0 &  0 \\
	0 & 0 & 0 & 0
\end{pmatrix}
+
\begin{pmatrix}
0 & 0 & 0 & 0 \\
0 & 0 & 0 & 0 \\
0 & 0 & 0 & 0\\
0 & 0 & 0 & K^{-1}
\end{pmatrix}
+
\begin{pmatrix}
0 & - \Dv & \alpha \grd & 0 \\
-\Grd & 0 & 0 & 0 \\ 
\alpha \dv & 0 & 0 & \dv\\ 
0 & 0 & \grd & 0
\end{pmatrix}
\right)
\begin{pmatrix}
	v\\\sigma\\p\\ q
\end{pmatrix}
= 
\begin{pmatrix}
	\rho f \\0\\g\\ 0
\end{pmatrix}\,.
\end{equation}
To further simplify the spatial differential operator in \eqref{Eq:HPS_20}, we introduce the total flux variable 
\begin{equation}
\label{Eq:HPS_21}
\olq := q + \alpha v
\end{equation}
and, then recast \eqref{Eq:HPS_20} as the evolutionary problem 
\begin{equation*}
	\label{Eq:HPS_22}
	\left(\partial_t 
	\begin{pmatrix}
		\rho & 0 & 0 &0\\
		0 & S & 0 & 0 \\
		0 & 0 & c_0 &  0 \\
		0 & 0 & 0 & 0
	\end{pmatrix}
	+
	\begin{pmatrix}
		-\alpha^2 K^{-1} & 0 & 0 & -\alpha K^{-1} \\
		0 & 0 & 0 & 0 \\
		0 & 0 & 0 & 0\\
	    -\alpha K^{-1} & 0 & 0 & K^{-1}
	\end{pmatrix}
	+
	\begin{pmatrix}
		0 & - \Dv & 0 & 0 \\
		-\Grd & 0 & 0 & 0 \\ 
		0 & 0 & 0 & \dv\\ 
		0 & 0 & \grd & 0
	\end{pmatrix}
	\right)
	\begin{pmatrix}
		v\\\sigma\\p\\ \olq
	\end{pmatrix}
	= 
	\begin{pmatrix}
		\rho f \\0\\g\\ 0
	\end{pmatrix}\,.
\end{equation*}
Finally, we put 
\begin{equation}
\label{Eq:HPS_30}
U := (v,\sigma,p,\olq)^\top \quad \text{and} \quad F := (\rho f,0, g, 0)^\top \,.
\end{equation}
We define the operators 
\begin{equation}
\label{Eq:HPS_31}
\hspace*{-1ex}
	M_0 := \begin{pmatrix}
	\rho & 0 & 0 &0\\
	0 & S & 0 & 0 \\
	0 & 0 & c_0 &  0 \\
	0 & 0 & 0 & 0
\end{pmatrix}\,, \;
M_1 := 	\begin{pmatrix}
	-\alpha^2 K^{-1} & 0 & 0 & -\alpha K^{-1}  \\
	0 & 0 & 0 & 0 \\
	0 & 0 & 0 & 0\\
	-\alpha K^{-1} & 0 & 0 & K^{-1}
\end{pmatrix}\,, \;
A:= \begin{pmatrix}
	0 & -\Dv & 0 & 0 \\
	-\Grd & 0 & 0 & 0 \\ 
	0 & 0 & 0 & \dv\\ 
	0 & 0 & \grd & 0
\end{pmatrix}\,.
\end{equation}
Then we obtain the following evolutionary problem.

\begin{prob}[Evolutionary problem]
\label{Prob:EP}
Let $H$ denote the product space 
\begin{equation}
\label{Eq:DefH}
H := L^2(\Omega)^d \times L^2(\Omega)^{d\times d}_{\operatorname{sym}} \times L^2(\Omega)\times L^2(\Omega)^d\,,
\end{equation}
equipped with the $L^2$ inner product of $L^2(\Omega)^{(d+1)^2}$. Let $M_0,M_1:H\rightarrow H$ and $A:D(A)\subset H\rightarrow H$, with 
\begin{equation}
\label{Eq:DefDA}
	D(A) := H^1_0(\Omega)^d \times D(\Dv) \times H^1_0(\Omega) \times D(\dv)\,,
\end{equation}
be defined by \eqref{Eq:HPS_31}. For given $F\in H_\nu(\R;H)$ according to \eqref{Eq:HPS_30}, find $U\in H_{\nu}(\R;H)$ such that 
\begin{equation}
\label{Eq:HPS_32}
(\partial_t M_0 + M_1 + A)U = F \,,
\end{equation}
where $U$ is defined by \eqref{Eq:HPS_30} along with \eqref{Def:NU}.
\end{prob}

Well-posedness of \eqref{Eq:HPS_32} is ensured by the following abstract result;  cf.\ \cite{P09} and \cite[Thm.~6.2.1]{STW22}.

\begin{thm}[Well-posedness]
\label{Thm:WP}
Let $H$ denote a real Hilbert space. Let $M_0,M_1:H\rightarrow H$ be bounded linear selfadjoint operators and $A:D(A)\subset H\rightarrow H$ skew-selfadjoint. Moreover, suppose that there exists some $\nu_0 >0$ such that	
\begin{equation}
\label{Eq:WP0}
\exists \, \gamma>0 \; \forall \nu \geq \nu_0\,, \; x\in H : \langle (\nu M_0+M_1)x,x\rangle_H \geq \gamma \langle  	x,x\rangle_H\,.
\end{equation}
Then, for each $\nu \geq \nu_0$ and each $F\in H_\nu(\R;H)$ there exist a unique solution $U\in H_{\nu}(\R;H)$ such that 
\begin{equation}
\label{Eq:WP1}
\overline{(\partial_t M_0 + M_1 + A)} U = F\,,
\end{equation}
where the closure is taken in $H_\nu(\R;H)$. Moreover, there holds the stability estimate 
\begin{equation*}
\label{Eq:WP2}
\| U \|_\nu \leq \frac{1}{\gamma} \| F \|_\nu \,. 
\end{equation*}
If $F\in H_\nu^s(\R;H)$ for some $s\in \N$, then the inclusion $U\in H_\nu^s(\R;H)$ is satisfied and the evolutionary equation is solved literally, such that
\begin{equation*}
	\label{Eq:WP3}
	{(\partial_t M_0 + M_1 + A)} U = F\,.
\end{equation*}
\end{thm}

\begin{cor}[Well-posedness of Problem \ref{Prob:EP}]
Problem  \ref{Prob:EP} is well-posed. In particular, there exists a unique solution $U\in H_{\nu}(\R;H)$ in the sense of \eqref{Eq:WP1}. 
\end{cor}

\begin{mproof}
For Problem~\ref{Prob:EP}, the assumptions of Thm.~\ref{Thm:WP} are fulfilled due to the conditions that the constants $\rho$, $\alpha$ and $c_0$ in \eqref{Eq:HPS_31} are strictly positive and the compliance tensor $S$ and matrix $K$ are symmetric and positive definite. The skew-selfadjointness of $A$ directly follows from \eqref{Def:div} and  \eqref{Def:Div}, respectively. Therefore, Thm.~\ref{Thm:WP} proves the assertion of this corollary. 
\end{mproof}

We note the following. 
\begin{rem}
\label{Rem:CP}
	
\begin{itemize}
\item By some version of the Sobolev embedding theorem (cf.~\cite[Lem.\ 3.1.59]{PMcG11}) there holds that 
\begin{equation}
\label{Eq:SET}
H^1_\nu(\R;H) \hookrightarrow C_\nu (\R;H)\,,
\end{equation}
where 
\begin{equation*}
C_\nu (\R;H) := \left\{f:\R \rightarrow H : f\; \text{is continuous}\,, \; \sup_{t\in \R} |f(t)|\operatorname{e}^{-\nu t}< \infty \right\}\,.	
\end{equation*}

\item For $F\in H^1_\nu(\R;H)$ there holds that $U\in H^1_\nu(\R;H)$ and, consequently, that  
\begin{equation}
\label{Eq:EmbDA}
AU = F - \partial_t M_0 U - M_1U \in H_\nu(\R;H)\,.
\end{equation}
Therefore, we have that $U(t)\in D(A)$ for almost every $t\in \R$. Moreover, for $F\in  H^2_\nu(\R;H)$ it follows that $U\in  H^2_\nu(\R;H)$ and, taking the time derivative of the equation in \eqref{Eq:EmbDA}, that $U\in  H^1_\nu(\R;D(A))$. By the embedding result \eqref{Eq:SET}, it then follows that $U\in C_\nu(\R;D(A))$. 

\item The condition \eqref{Eq:WP0} of positive definiteness is assumed to hold uniformly in $\nu \geq \nu_0$. This ensures the causality of the solution operator $S_\nu := \overline{(\partial_t M_0 + M_1 + A)} ^{-1}$ to \eqref{Eq:WP1}; cf.~\cite[Thm.~6.2.1]{STW22}. Further, it holds that $S_\nu F = S_\eta F$ for $\nu,\eta \geq \nu_0$ and $F\in H_\nu(\R;H)\cap H_\eta (\R;H)$.

\item Initial value problems for \eqref{Eq:HPS_32} are studied by a generalization of the solution theory to certain distributional right-hand sides; cf.~\cite[Thm.~9.4.3]{STW22} and \cite[Thm.~6.2.9]{PMcG11}. Let $F\in H_\nu (\R;H)$ be supported on $[t_0,\infty)$ for some $t_0\in \R$ and $U_0\in D(A)$ be given. Then, the evolutionary equation 
\begin{equation}
\label{Eq:EPpIC}
(\partial_t M_0 + M_1 + A) U = F + \delta_{t_0} M_0 U_0\,,	
\end{equation} 
has a unique solution $U\in H_\nu^{-1}(\R;H)$ satisfying $MU(t_0^+)=M_0 U_0$ in $H^{-1}(D(A))$, where $\delta_{t_0}$ denotes the delta distribution at $t=t_0$. For $t\in (0,\infty)$, the evolutionary equation $(\partial_t M_0 + M_1 + A) U = F$ is satisfied in the sense of distributions for test functions $\varphi \in H^1_\nu (\R;H)\cap H_\nu(\R;D(A))$ supported on $[t_0,\infty)$. The distribution on the right-hand side of \eqref{Eq:EPpIC} can still be avoided by reformulating the initial value problem for \eqref{Eq:WP3} into the evolutionary equation 
\begin{equation*}
(\partial_t M_0 + M_1 + A) W = \partial_t^{-1} F + H_{t_0}U_0
\end{equation*}
for $W:=\partial_t^{-1}U$, where $H_{t_0}$ denotes the Heaviside function with jump in $t=t_0$; cf.\ \cite[Cor.\ 1.1]{F21} and \cite[p.\ 446]{PMcG11}. Then, it follows that
\begin{equation*}
	U = \partial_t \overline{(\partial_t M_0 + M_1 + A)}^{\;-1} \partial_t^{-1} (F+\delta_{t_0} U_0)\,.
\end{equation*}

\end{itemize}
\end{rem}

\section{Discontinuous Galerkin discretization and well-posedness}
\label{Sec:TimeDisc}

Here we derive a family of fully discrete schemes for Problem~\ref{Prob:EP}. Space and time discretization are based on discontinuous Galerkin approaches. Well-posedness of the discrete problem is shown. We assume that the weight $\nu$ in \eqref{Eq:Hnu} is chosen such that the assumptions of Thm.~\ref{Thm:WP} are satisfied.

\subsection{Notation and auxiliaries}

For the time discretization, we decompose $I=(0,T]$ into $N$ subintervals $I_n=(t_{n-1},t_n]$, for $n=1,\ldots,N$, where $0=t_0<t_1< \cdots < t_{N-1} < t_N = T$ such that $I=\bigcup_{n=1}^N I_n$. We put $\tau := \max_{n=1,\ldots, N} \tau_n$ with $\tau_n = t_n-t_{n-1}$. Further, the set $\mathcal{M}_\tau := \{I_1,\ldots, I_N\}$ is called the time mesh. For any $k\in \N_0$ and some Banach space $B$, we let 
\begin{equation}
	\label{Def:Pk}
	\mathbb P_k(I_n;B) := \bigg\{w_\tau \,: \,  I_n \to B \,, \mid w_\tau(t) = \sum_{j=0}^k 
	W^j t^j \; \forall t\in I_n\,, \; W^j \in B\; \forall j \bigg\}\,
\end{equation}
denote the space of $B$-valued polynomials of degree at most $k$ defined on $I_n$. For a Hilbert space $H$, the space $\mathbb P_k(I_n;H)$, equipped with the exponentially weighted in time inner product
\begin{equation}
\label{Def:IPPk}
\langle v_\tau, w_\tau\rangle_{\nu,n} : = \int_{t_{n-1} }^{t_n} \langle v_\tau(t), w_\tau(t)\rangle_H \operatorname{e}^{-2\nu (t-t_{n-1})} \ud t \,,
\end{equation}
is a Hilbert space.  The semidiscretization in time of \eqref{Eq:HPS_32} by Galerkin methods is done in
\begin{equation}
\label{Eq:DefYk} 
Y_{\tau,\nu}^k (B) := \left\{w_\tau \in H_\nu (0,T;B) \mid  w_\tau{}_{|I_n} \in \mathbb P_k(I_n;B)\; \forall I_n\in \mathcal{M}_\tau\,, \; w_\tau(0) \in B \right\}\,.
\end{equation}
For any function $w: \overline I\to B$ that is piecewise sufficiently smooth with respect to the time mesh $\mathcal{M}_{\tau}$, for instance for $w\in Y^k_\tau (B)$, we define the right-hand sided and left-hand sided limit at a mesh point $t_n$ by
\begin{equation}
	\label{Eq:DefLim}
	w^+(t_n) := \lim_{t\to t_n+0} w(t) ,\quad\text{for}\; n<N\,,
	\qquad\text{and}\qquad
	w^-(t_n) := \lim_{t\to t_n-0} w(t) ,\quad\text{for}\; n>0\,.
\end{equation}
For the error analysis, we further need the space 
\begin{equation}
	\label{Eq:DefXk} 
	X_{\tau,\nu}^k (B) := \left\{w_\tau \in C_\nu ([0,T];B) \mid  w_\tau{}_{|I_n} \in \mathbb P_k(I_n;B)\; \forall I_n\in \mathcal{M}_\tau\right\}\,.
\end{equation}

In the discrete scheme, a quadrature formula is applied for the evaluation of the time integrals. For the discontinuous in time finite element method, a natural choice is to consider the $(k+1)$-point right-sided Gau{ss}--Radau quadrature formula on each time interval $I_n=(t_{n-1},t_n]$. Here, we use a modification of the standard right-sided Gau{ss}--Radau quadrature formula that is defined by 
\begin{equation}
	\label{Eq:GRF}
	Q_{n,\nu}(w) := \frac{\tau_n}{2}\sum_{\mu=0}^{k} \hat \omega_\mu w(t_{n,\mu}) \approx  \int_{I_n} \operatorname{e}^{-2\nu (t-t_{n-1})}  w(t) \ud t\,,  
\end{equation}
where $t_{n,\mu}=T_n(\hat t_{\mu})$, for $ \mu = 0,\ldots ,k$, are the quadrature points on $I_n$ and $\hat  \omega_\mu$ the corresponding weights. Here, $T_n(\hat t\,):=(t_{n-1}+t_n)/2 + (\tau_n/2)\hat t$ is the affine transformation from the reference interval $\hat I = (-1,1]$ to $I_n$ and $\hat t_{\mu}$, for $\mu = 0,\ldots,k$, are the quadrature points of the weighted Gau{ss}--Radau formula on $\hat I$ (cf.~\cite{PTVF07}), such that for all polynomials $p\in \mathbb P_{2k}(\hat I;\R)$ there holds that 
\begin{equation*}
	\label{Eq:GRFU}
	 \int_{\hat I}  \operatorname{e}^{-\nu \tau_n (\hat t+1)}  p(\hat t\,) \ud \hat t = \sum_{\mu=0}^{k} \hat \omega_\mu p(\hat t_{\mu}) \,.
\end{equation*}
Then, for polynomials $p\in \mathbb P_{2k}(I_n;\R)$ we have that 
\begin{equation}
\label{Eq:ExactQuad}
	Q_{n,\nu}(p) = 	 \int_{I_n} \operatorname{e}^{-2\nu (t-t_{n-1})}  p(t) \ud t \,.
\end{equation}
Finally, we introduce the time-mesh dependent quantities 
\begin{subequations}
\begin{alignat}{4}
Q_n[w]_\nu &  := Q_{n,\nu} (w)\,, & Q_n[v,w]_\nu  & := Q_{n,\nu} (\langle v,w\rangle_H )\,,\\[1ex]
\label{Eq:TdAV}
| w | _{\tau,\nu,n}^2 & := Q_n[w]_\nu\,, & | w |_{\tau,\nu}^2 & := \sum_{n=1}^N Q_n[w]_\nu \operatorname{e}^{-2\nu t_{n-1}}\,, \\[1ex]
\label{Eq:TdN}
\| w \| _{\tau,\nu,n}^2 & := Q_n[w,w]_\nu\,, & \| w \|_{\tau,\nu}^2 & := \sum_{n=1}^N Q_n[w,w]_\nu \operatorname{e}^{-2\nu t_{n-1}}\,,
\end{alignat}
\end{subequations}
where the nonnegativity of $w$ is tacitly assumed in the definition of $| w | _{\tau,\nu,n}^2$. This will be satisfied below. 

For the nodes $t_{n,\mu}\in (t_{n-1},t_n]$, for $\mu=0,\ldots,k$ and $n=1,\ldots,N$, of the weighted Gauss--Radau formula \eqref{Eq:GRF}, we define the global Lagrange interpolation operator $I_\tau: C([0,T];B)\rightarrow Y_\tau^k(B)$ by
\begin{equation}
\label{Def:LIO}
I_\tau f(0) = f(0) \,, \qquad I_\tau f(t_{n,\mu}) = f(t_{n,\mu})\,, \quad \mu=0,\ldots,k\,, \; n=1,\ldots,N\,.
\end{equation}
For the Lagrange interpolation \eqref{Def:LIO}, on each $I_n$ there holds that (cf.~\cite[Thm.~1]{H91})
\begin{equation}
	\label{Eq:ErrLI}
	\| w - I_\tau w \|_{C(\overline I_n;B)} \leq c \tau_n^{k+1} \| \partial_t^{k+1} w\|_{C(\overline I_n;B)}\,.
\end{equation}
Moreover, we need the Lagrange interpolation operator $\hItau: C([0,T];B)\rightarrow X_\tau^{k+1}(B)$ with respect to the Gauss--Radau quadrature points $t_{n,\mu}$, for $\mu=0,\ldots,k$, and $t_{n-1}$, for $n=1,\ldots,N$, that is defined by 
\begin{equation}
	\label{Def:LIOc}
	\hItau f(t_{n,\mu}) = f(t_{n,\mu})\,, \quad \mu=0,\ldots,k\,,  \quad \text{and} \quad \hItau f(t_{n-1}) = f(t_{n-1}) \quad \text{for} \;\; n=1,\ldots,N\,.
\end{equation}
Then, for $\hItau$ there holds that (cf.~\cite[Thm.~2]{H91})
\begin{equation}
	\label{Eq:ErrLI0c}
	\| \partial_t^{s} (w - \hItau w) \|_{C(\overline I_n;B)} \leq c \tau_n^{k+2-s} \| \partial_t^{k+2} w\|_{C(\overline I_n;B)}\,, \quad  \text{for }\; s\in \{0,1\}\,.
\end{equation}

For the space discretization, let the mesh $\mathcal{T}_h=\{K\}$ denote a decomposition of the polyhedron $\Omega$ into  quadrilateral or hexahedral elements $K$ with meshsize $h=\max\{h_K \; : \; K\in \mathcal T_h\}$ for  $h_K:=\operatorname{diam}(K)$. The mesh is assumed to be conforming (matching) and shape-regular; cf., e.g., \cite{R08}. The assumptions about $\mathcal{T}_h$ are sufficient to derive inverse and trace inequalities; cf.~\cite[Chap.~1]{PE12}. Further, optimal polynomial approximation properties in the sense of \cite[Def.~1.55]{PE12} are satisfied; cf., e.g., \cite[Thm.~2.6]{R08}. Simplicial triangulations can be considered analogously. For more general mesh concepts in the context of discontinuous Galerkin methods we refer to \cite[Subsec.~1.4]{PE12} or \cite[Subsec.~2.3.2]{DF15}. For any $K\in \mathcal T_h$ we denote by $n_K$ the outward unit normal to the faces (egdes for $d=2$) of $K$. Further, we let $\mathcal E_h$ be the union of the boundaries of all elements of $\mathcal T_h$. Let $\mathcal E_h^i = \mathcal E_h \backslash \partial \Omega$ be the set of interior faces (edges if $d=2$) and $\mathcal E_h^\partial =\mathcal E_h\backslash \mathcal E_h^i$ denote the union of all boundary faces.  

For any $r\in \N$, the discrete space of continuous and piecewise polynomial functions is denoted as 
\begin{equation}
	\label{Eq:DefXhr}
	X_h^r :=  \{ w_h \in C(\overline \Omega) \mid w_h{}_{|K} \in W_r(K)\; \; \forall K\in \mathcal T_h \}\cap H^1_0(\Omega)\,,
\end{equation} 
where the local space $W_r(K)$ is defined by mapped versions of $\mathbb Q_{r}$; cf.\ \cite[Subsec.\ 3.2]{QV08}. For any $r\in \N_0$, we denote the space of broken polynomials by 
\begin{equation}
	\label{Eq:DefYhr}
	Y_h^r :=  \{ w_h \in L^2(\Omega) \mid w_h{}_{|K} \in W_r(K)\; \; \forall K\in \mathcal T_h \}\,.
\end{equation}
For the spatial approximation of Problem~\ref{Prob:EP} we consider using 
\begin{equation}
	\label{Eq:DefHh}
	\Hhnc\in \{H_h^{\text{hy}},H_h^{\text{dg}}\}\,, \qquad \text{with }\; \Hhnc\subset H\,,
\end{equation}
where the finite element product spaces  $H_h^{\text{hy}}$ and $H_h^{\text{dg}}$ are given by  
\begin{subequations}
	\label{Eq:DefHhhydg}
	\begin{alignat}{2}
		\label{Eq:DefHhhy}
		H_h^{\text{hy}}  & := (X_h^r)^d \times \big((Y_h^r)^{d\times d}\cap L^2(\Omega)^{d\times d}_{\operatorname{sym}}\big) \times X_h^r \times (Y_h^r)^d\,, \\[1ex]
		\label{Eq:DefHhdg}
		H_h^{\text{dg}}  & := (Y_h^r)^d \times \big((Y_h^r)^{d\times d}\cap L^2(\Omega)^{d\times d}_{\operatorname{sym}}\big) \times Y_h^r \times (Y_h^r)^d\,. 
	\end{alignat}
\end{subequations}
Discretizations of Problem~\ref{Prob:EP} in either spaces, $H_h^{\text{hy}}$ and $H_h^{\text{dg}}$, are studied simultaneously. The reason for considering also the hybrid space $H_h^{\text{hy}}$ is that continuous and $H^1_0(\Omega)$-conforming finite element methods lead to lower computational cost than discontinuous ones. $H(\operatorname{div};\Omega)$-conforming approximations in the framework of Picard's theory have been studied in \cite{FTW19} for scalar-valued problems of changing type. These families of schemes can be applied analoguously to the approximation of $\sigma$ and $\overline q$ in Problem~\ref{Prob:EP}. Since discontinuous Galerkin methods offer high flexibilty combined with  implementational advantages, DG methods are attractive and studied here.  

In Subsec.~\ref{Subsec:ErrEst} we need the $L^2$-orthogonal projection of functions $w\in H$ onto the broken polynomial space $H_h^{\text{dg}}$ of \eqref{Eq:DefHhdg} that is very simple, even on more general meshes than studied here. For the $L^2$-orthogonal projection $\Pi_h: H \to H_h^{\text{dg}}$,  $v\in H$ and $\Pi_h v \in H_h^{\text{dg}}$, with 
\begin{equation}
	\label{Def:Pih}
	\langle \Pi_h v, w_h \rangle_H = \langle v, w_h \rangle_H\,, \quad \text{for all} \; w_h \in H_h^{\text{dg}}\,, 
\end{equation}
there holds for all $s\in \{0,\ldots,r+1\}$ and all $w\in H^s(K)$ that 
\begin{equation}
	\label{Eq:IOPih} 
	| w - \Pi_h w|_{H^m (K)} \leq C_{\operatorname{app}} h_K^{s-m} |w|_{H^s(K)}\,, \quad \text{for }\; m\in \{0,\ldots,s\}\,,
\end{equation}
where $C_{\operatorname{app}}$ is independent of both $K$ and $h$; cf.\ \cite[Lem.~1.58]{PE12}, \cite[Thm.\ 2.6]{R08}.  In \eqref{Eq:IOPih}, we denote by $|\cdot |_{H^m (K)}$ the seminorm of the Sobolev space $H^m (K)$. Also, the $L^2$-orthgonal projection satisfies that 
\begin{subequations}
	\label{Eq:IOPih2} 
	\begin{alignat}{2}
		\label{Eq:IOPih21} 
		\| w - \Pi_h w\|_{L^2(e)} & \leq C'_{\operatorname{app}} h_K^{s-1/2} |w|_{H^s(K)}\,, \quad \text{for }\; s\geq 1\,,\\[1ex]
		\label{Eq:IOPih22} 
		\| \nabla (w - \Pi_h w)_{|K}\cdot n_e \|_{L^2(e)} & \leq C''_{\operatorname{app}} h_K^{s-3/2} |w|_{H^s(K)}\,, \quad \text{for }\; s\geq 2\,,
	\end{alignat}	
\end{subequations}
where $C'_{\operatorname{app}}$ and $C''_{\operatorname{app}}$ are independent of both $K$ and $h$; cf.\ \cite[Lem.~1.59]{PE12}.

 \subsection{Gradient and divergence on broken function spaces}

To define our discontinuous Galerkin discretization schemes we need to recall some general concepts for the definition of the gradient and divergence operator on broken function spaces with respect to the triangulation $\mathcal T_h$. For further details and concepts of broken function spaces we refer to, e.g., \cite{PE12}. On the triangulation $\mathcal T_h$, let $Y_h=Y_h(\mathcal T_h)$ and $Z_h=Z_h(\mathcal T_h)$ denote piecewise (broken) spaces of scalar- and vector-valued functions, respectively. On the set of (inner and outer) boundaries $\mathcal E_h$, let $\widehat Y_h=\widehat Y_h(\mathcal E_h)$ and $\widehat Z_h = \widehat Z_h(\mathcal E_h)$ be piecewise (broken) spaces of scalar- and vector-valued functions on $\mathcal E_h$, respectively. We put $\widetilde Y_h := Y_h\times \widehat Y_h$ and $\widetilde Z_h := Z_h\times \widehat Z_h$. We denote the dual spaces of $\widetilde Y_h$ and $\widetilde Z_h$ by $\widetilde Y_h^\ast$ and $\widetilde Z_h^\ast$. In these spaces we define the following derivatives of the discontinuous Galerkin method; cf.~\cite{HWWX19}. 

\begin{defi}[DG derivatives]
\label{Def:dgOp}
Let $\widetilde y_h = (y_h,\widehat y_h) \in \widetilde Y_h $ and $\widetilde z_h = (z_h,\widehat z_h)\in \widetilde Z_h$. Then the DG-gradient $\dggrd: Y_h \rightarrow \widetilde Z_h^\ast$ and  the DG-divergence $\dgdiv: Z_h \rightarrow \widetilde Y_h^\ast$ are defined by  
\begin{subequations}
\label{Eq:dgOp}
\begin{alignat}{2}
\label{Eq:dggrd}
\langle  \dggrd y_h , \widetilde z_h\rangle & := \langle \operatorname{grad}_h y_h, z_h \rangle - \sum_{K\in \mathcal T_h} \langle y_h, \widehat z_h \cdot n_K \rangle_{\partial K} \,, & \qquad \forall y_h \in Y_h\,, \quad \forall  \widetilde  z_h \in \widetilde Z_h\,, \\[1ex]
\label{Eq:dgdiv}
\langle  \dgdiv z_h , \widetilde y_h\rangle & := \langle \operatorname{div}_h z_h, y_h \rangle - \sum_{K\in \mathcal T_h} \langle z_h \cdot n_K , \widehat y_h \rangle_{\partial K}  \,, & \qquad \forall z_h \in Z_h\,, \quad  \forall \widetilde y_h \in \widetilde Y_h\,.
\end{alignat}
\end{subequations}
\end{defi}

Here, $\langle \cdot ,  \cdot \rangle_S$ denotes the inner product of $L^2(S)$, where we drop the index $S$ if $S= \Omega $. Further,  $\operatorname{grad}_h $ and $\operatorname{div}_h$ are the broken gradient and divergence, respectively; cf.~\cite[Subsec.~1.2.5 and 1.2.6]{PE12}. In what follows, we drop the index $h$ in the broken operators, when this operation appears inside an integral over a fixed mesh element $K\in \mathcal T_h$. We recall that on the usual Sobolev spaces the broken gradient coincides with the distribution gradient; cf.~\cite[Lem.~1.22]{PE12}. The same applies to the broken divergence; cf.~\cite[Subsec.~1.2.6]{PE12}. The dual operators of $\dggrd $ and $\dgdiv$ are denoted by $\dggrd ^\ast: \widetilde Z_h \rightarrow Y_h^\ast$ and $\dgdiv: \widetilde Y_h \rightarrow Z_h^\ast$. Then, there holds that
 \begin{subequations}
 	\label{Eq:ddgOp}
 	\begin{alignat}{2}
 		\label{Eq:ddggrd}
 		\langle  \dggrd^\ast \widetilde z_h , y _h\rangle & = \langle \widetilde z_h , \dggrd y_h \rangle\,, & \qquad  \qquad \forall \widetilde z_h  \in  \widetilde Z_h  \,,  & \quad \forall y_h \in Y_h\,,   \\[1ex]
 		\label{Eq:ddgdiv}
 		\langle  \dgdiv^\ast \widetilde y_h, z_h\rangle & = \langle \widetilde y_h,  \dgdiv z_h \rangle \,, & \qquad  \forall   \widetilde y_h \in  \widetilde Y_h \,, & \quad \forall z_h \in Z_h  \,.
 	\end{alignat}
 \end{subequations}

The DG-derivatives $ \dggrd$ and $-\dgdiv$ are \textit{conditionally dual} with each other; cf.~\cite{HWWX19}. To demonstrate this link, we deduce from \eqref{Eq:dgOp} and \eqref{Eq:ddgOp} that 
\begin{subequations}
\label{Eq:DOPdgdiv}
\begin{alignat}{2}
\nonumber
\langle - \dgdiv^\ast \widetilde y_h , z_h \rangle & = - \langle \widetilde y_h, \dgdiv z_h \rangle \\
& = \langle \operatorname{grad}_h y_h, z_h  \rangle + \sum_{K\in \mathcal T_h} \langle \widehat y_h - y_h , z_h \cdot n_K \rangle_{\partial K}\,, \\[1ex]
\langle \dggrd y_h, \widetilde z_h \rangle & = \langle \operatorname{grad}_h y_h, z_h \rangle -   \sum_{K\in \mathcal T_h} \langle y_h , \widehat z_h \cdot n_K \rangle_{\partial K}
\end{alignat}
\end{subequations}
and
\begin{subequations}
	\label{Eq:DOPdggrd}	
	\begin{alignat}{2}
		\nonumber
		\langle \dggrd^\ast \widetilde z_h , y_h \rangle & = \langle \widetilde z_h, \dggrd y_h \rangle \\
		& = - \langle \operatorname{div}_h z_h, y_h  \rangle + \sum_{K\in \mathcal T_h} \langle (z_h - \widehat z_h)\cdot n_K , y_h \rangle_{\partial K}\,, \\[1ex]
		- \langle \dgdiv z_h, \widetilde y_h \rangle & = - \langle \operatorname{div}_h  z_h , y_h \rangle +   \sum_{K\in \mathcal T_h} \langle z_h \cdot n_K,  \widehat y_h \rangle_{\partial K}\,.
	\end{alignat}
\end{subequations}
The identities \eqref{Eq:DOPdgdiv} and \eqref{Eq:DOPdggrd} directly imply the following conditional duality between $\dggrd$ and $\dgdiv$ under the assumption that $\widehat y_h = 0$ on $\mathcal E_h^\partial$; cf.\ also \cite[Lem.~2.1]{HWWX19}.

\begin{lem}
\label{Lem:CondDual}
Suppose that $\widehat y_h = 0$ on $\mathcal E_h^\partial$. For the DG derivatives \eqref{Eq:dgOp} there holds the duality
\begin{equation}
\label{Eq:Discdual}	
\dgdiv = - (\dggrd)^\ast\,,
\end{equation}
if one of the following conditions is satisfied:
\begin{subequations}
\label{Eq:DiscdualCond}	
\begin{alignat}{2}
\label{Eq:DiscdualCond1}	
\text{i)}   \hspace*{2ex}\mbox{} & z_h\cdot n_K{}_{|\mathcal E_h} = \widehat z_h\cdot n_K\,; \hspace*{30ex}\mbox{}\\[1ex]
\label{Eq:DiscdualCond2}	
\text{ii)}  \hspace*{2ex}\mbox{} & y_h{}_{|\mathcal E_h} = \widehat y_h\,;\\[1ex]
\label{Eq:DiscdualCond3}	
\text{iii)} \hspace*{2ex}\mbox{} & \widehat z_h\cdot n_K = \frac{1}{2} (z_h^+ + z_h^-) \cdot n_K \quad \text{and} \quad \widehat y_h = \frac{1}{2} (y_h^+ + y_h^-)\,.
\end{alignat}
\end{subequations}	
\end{lem}
In \eqref{Eq:DiscdualCond3}, we let $y_h^{\pm} := y_h{}_{|\partial K^\pm}$ and $z_h^{\pm} := z_h{}_{|\partial K^\pm}$ for two adjacent elements $K^+$ and $K^-$ with common face $e\in  \mathcal E_h^i$ and outer unit normal vector $n^\pm$ to $\partial K^\pm$. By $n_K$ we denote the normal vector assigned to $e\in \mathcal E_h$, where $n_K$ is the outer normal vector for $e\in \mathcal E_h^\partial$.  We note that $\dggrd y_h = \operatorname{grad} y_h$ for $y_h \in H^1_0(\Omega)$ and $\dgdiv z_h = \operatorname{div} z_h$ for $z_h \in H(\operatorname{div},\Omega)$, since the second terms on the right-hand side of \eqref{Eq:dgOp} yield that 
\begin{align*}
\sum_{K\in \mathcal T_h} \langle y_h, \widehat z_h \cdot n_K \rangle_{\partial K} & = \sum_{e\in \mathcal E_h^i} \langle y_h, \widehat z_h \cdot ( n^+ + n^-) \rangle_{e} = 0 \,,\\[1ex]
\sum_{K\in \mathcal T_h} \langle z_h \cdot n_K , \widehat y_h \rangle_{\partial K} & = \sum_{e\in \mathcal E_h^i} \langle z_h^+ \cdot n^+ + z_h^- \cdot  n^-,  \widehat y_h \rangle_{e} = 0\,.
\end{align*}
The matrix- and vector-valued operators $\Grd$ and $\Dv$, introduced in \eqref{Def:Grad0} and \eqref{Def:Div} respectively, are defined on broken function spaces  similarly to \eqref{Eq:dgOp}.

\subsection{Gradient and divergence on broken polynomial spaces and operator~$\vec{A_h}$}
\label{Subsec:NconfFEM}

Now, we specify the broken spaces $\widehat Y_h$ and $\widehat Z_h$ of Def.~\ref{Def:dgOp} for the finite element spaces \eqref{Eq:DefHh} and \eqref{Eq:DefHhhydg} that we consider for the spatial approximation of Problem~\ref{Prob:EP}. In light of Lem.~\ref{Lem:CondDual} we put
\begin{subequations}
\label{Eq:DefTraceFEM}
\begin{alignat}{7}
\label{Eq:DefTraceFEM_1}
	\widehat y_h & :=   \frac{1}{2} (y_h^+ + y_h^-) & & \quad \text{and} &  \quad  \widehat z_h\cdot n_e & := \frac{1}{2} (z_h^+ + z_h^-) \cdot n_e \,,& \quad \text{for} \;\; e\in \mathcal E_h^i\,, \\[1ex]
\label{Eq:DefTraceFEM_2}
	\widehat y_h & := y_h  \,,  && \quad \text{and} &  \quad \widehat  z_h\cdot n_e & :=  z_h\cdot n_e & \quad \text{for} \;\; e\in \mathcal E_h^\partial\,.
\end{alignat}
\end{subequations}
For \eqref{Eq:DefTraceFEM}, the definitions of the DG gradient operator $\dggrd: Y_h^r \to {((Y_h^r)^{d})}^\ast$ in \eqref{Eq:dggrd} and the DG divergence operator $\dgdiv: (Y_h^r)^{d}\to ({Y_h^r})^\ast$ in \eqref{Eq:dgdiv} of Def.~\ref{Def:dgOp}  then read as follows. 


\begin{defi}[DG derivatives for single-valued functions]
\label{Def:divgrdDG}
The DG gradient operator $\dggrd: Y_h^r \to {((Y_h^r)^{d})}^\ast$ and the DG divergence operator $\dgdiv: (Y_h^r)^{d}\to ({Y_h^r})^\ast$ are defined by 
\begin{subequations}
	\label{Eq:dgOpHh}
	\begin{alignat}{2}
		\label{Eq:dggrdHh}
		\langle  \dggrd y_h , z_h\rangle & := \langle \operatorname{grad}_h y_h, z_h \rangle - \sum_{e\in \mathcal E_h^i} \langle \lsem y_h \rsem, \ldblbrace z_h\rdblbrace  \cdot n_e \rangle_{e} - \sum_{e\in \mathcal E_h^\partial} \langle y_h ,  z_h   \cdot n_e \rangle_{e} \,, \\[1ex]
		\label{Eq:dgdivHhl}
		\langle  \dgdiv z_h , y_h\rangle & := \langle \operatorname{div}_h z_h, y_h \rangle - \sum_{e\in \mathcal E_h^i} \langle \lsem z_h \rsem \cdot n_e , \ldblbrace y_h \rdblbrace \rangle_{e}  - \sum_{e\in \mathcal E_h^\partial} \langle z_h  \cdot n_e , y_h \rangle_{e} 
	\end{alignat}
\end{subequations}
for all $y_h \in Y_h^{r}$ and $z_h \in (Y_h^r)^d$,  where standard notation (cf.~\cite{PE12}) is used for the averages and jumps 
\begin{equation*}
\ldblbrace w \rdblbrace_e :=  \frac{1}{2} (w^{+}+ w^{-})\quad \text{and} \quad  \lsem w \rsem_e : =  w^{+} -  w^{-} \,.
\end{equation*}
\end{defi}

On the usual Sobolev spaces the DG gradient and DG divergence of \eqref{Eq:dgOpHh} coincide with the distribution gradient and divergence, respectively, since for functions of $H_0^1(\Omega)$ the jump terms $\lsem y_h \rsem$ on $e\in \mathcal E_h^i$ and the traces on the boundary faces $e\in \mathcal E_h^\partial$ vanish in \eqref{Eq:dgOpHh}; cf.~\cite[Lem.~1.22 and 1.24]{PE12}.  Similarly, for functions of $H(\operatorname{div};\Omega)$ the jumps $\lsem z_h \rsem \cdot n_e$ for $e\in \mathcal E_h^i$ vanish as well; cf.~\cite[Lem.~1.22 and 1.24]{PE12}.  The assumption of Lem.~\ref{Lem:CondDual} that $\widehat y_h = 0$ for $e\in \mathcal E_h^\partial$ is not fulfilled for the DG space \eqref{Eq:DefHh} with \eqref{Eq:DefHhdg}. This leads to a perturbation of the skew-selfadjointness of $\dggrd$ and $\dgdiv$  shown now. 

\begin{lem}[DG skew-selfadjointness]
\label{Lem:CondDualIBV}
For all $y_h \in Y_h^{r}$ and $z_h \in (Y_h^r)^d$ there holds that 
\begin{subequations}
\label{Eq:CondDualPtb0}
\begin{alignat}{2}
\label{Eq:CondDualPtb1}
\langle  \dggrd y_h , z_h\rangle  + \sum_{e\in \mathcal E_h^\partial} \langle y_h ,  z_h   \cdot n_e \rangle_{e} & = \langle  - \dgdiv^\ast y_h , z_h\rangle \,,\\[1ex]
\label{Eq:CondDualPtb2}
\langle  \dgdiv z_h , y_h\rangle + \sum_{e\in \mathcal E_h^\partial} \langle y_h , z_h \cdot n_e \rangle_e & = \langle - \dggrd^\ast z_h, y_h \rangle
\end{alignat}
\end{subequations}
and 
\begin{equation}
\label{Eq:CondDualPtb3}
\langle  \dggrd y_h , z_h\rangle + \langle  \dgdiv z_h , y_h\rangle = - \sum_{e\in \mathcal E_h^\partial} \langle y_h , z_h \cdot n_e \rangle_e\,.
\end{equation}
\end{lem}	

\begin{mproof}
The identity \eqref{Eq:CondDualPtb0} is a direct consequence of \eqref{Eq:DOPdgdiv} and \eqref{Eq:DOPdggrd} along with the definition \eqref{Eq:DefTraceFEM} of the broken spaces $\widehat Y_h$ and $\widehat Z_h$ on $\mathcal E_h$. From \eqref{Eq:CondDualPtb2} we then get that 
\begin{align*}
\langle  \dggrd y_h , & z_h\rangle  + \langle  \dgdiv z_h , y_h\rangle = \langle  \dggrd^\ast z_h, y_h\rangle + \langle  \dgdiv z_h , y_h\rangle\\[1ex]
& = - \langle  \dgdiv z_h , y_h\rangle - \sum_{e\in \mathcal E_h^\partial} \langle y_h , z_h \cdot n_e \rangle_e  + \langle  \dgdiv z_h , y_h\rangle  = - \sum_{e\in \mathcal E_h^\partial} \langle y_h , z_h \cdot n_e \rangle_e\,.
\end{align*}
This proves \eqref{Eq:CondDualPtb3}.

\end{mproof}

DG derivatives $\dgGrd: (Y_h^r)^d \to \big( (Y_h^r)^{d\times d}\cap L^2(\Omega)^{d\times d}_{\operatorname{sym}}\big)^\ast$ and $\dgDiv: (Y_h^r)^{d\times d}\cap L^2(\Omega)^{d\times d}_{\operatorname{sym}}  \to \big((Y_h^r)^{d}\big)^\ast$ for multi-valued function are defined as follows.

\begin{defi}[DG derivatives for multi-valued functions]
For multi-valued functions, the DG gradient operator $\dgGrd: (Y_h^r)^d \to \big( (Y_h^r)^{d\times d}\cap L^2(\Omega)^{d\times d}_{\operatorname{sym}}\big)^\ast$ and DG divergence operator $\dgDiv: (Y_h^r)^{d\times d}\cap L^2(\Omega)^{d\times d}_{\operatorname{sym}}  \to \big((Y_h^r)^{d}\big)^\ast$ are given by 
\begin{subequations}
	\label{Eq:dgOpHhV}
	\begin{alignat}{2}
		\label{Eq:dggrdHhV}
		\langle  \dgGrd  y_h , z_h\rangle & := \langle \operatorname{Grad}_h y_h, z_h \rangle - \sum_{e\in \mathcal E_h^i} \langle \lsem y_h \rsem, \ldblbrace z_h\rdblbrace  \cdot n_e \rangle_{e} - \sum_{e\in \mathcal E_h^\partial} \langle y_h,  z_h \cdot n_e \rangle_{e}   \,, \\[1ex]
		\label{Eq:dgdivHhlV}
		\langle  \dgDiv z_h , y_h\rangle & := \langle \operatorname{Div}_h z_h, y_h \rangle - \sum_{e\in \mathcal E_h^i}  \langle \lsem z_h \rsem \cdot n_e , \ldblbrace y_h \rdblbrace \rangle_{e}  - \sum_{e\in \mathcal E_h^\partial} \langle z_h \cdot n_e , y_h \rangle_{e}
	\end{alignat}
\end{subequations}
for all $y_h \in (Y_h^r)^d$ and $z_h \in (Y_h^r)^{d\times d}\cap L^2(\Omega)^{d\times d}_{\operatorname{sym}} $. 
\end{defi}

In \eqref{Eq:dgOpHhV},  the operators $\operatorname{Grad}_h$ and $\operatorname{Div}_h $ are the broken symmetrized gradient and broken divergence that extend the distributional gradient in \eqref{Def:Grad0} and divergence in \eqref{Def:Div} to broken polynomial spaces; cf.\ \cite[Def.~1.21]{PE12}. On the usual Sobolev spaces the broken gradient $\dgGrd$ and divergence $\dgDiv$ coincide with the distributional symmetrized gradient and divergence of \eqref{Def:Grad0} and \eqref{Def:Div}, respectively. Similarly to Lem.~\ref{Lem:CondDualIBV}, for DG spaces there holds for $y_h \in (Y_h^r)^d$ and $z_h \in (Y_h^r)^{d\times d}\cap L^2(\Omega)^{d\times d}_{\operatorname{sym}} $ that
\begin{subequations}
	\label{Eq:CondDualPtb10}
	\begin{alignat}{2}
		\label{Eq:CondDualPtb11}
		\langle  \dgGrd y_h , z_h\rangle  + \sum_{e\in \mathcal E_h^\partial} \langle y_h ,  z_h   \cdot n_e \rangle_{e} & = \langle  - \dgDiv^\ast y_h , z_h\rangle \,,\\[1ex]
		\label{Eq:CondDualPtb12}
		\langle  \dgDiv z_h , y_h\rangle + \sum_{e\in \mathcal E_h^\partial} \langle y_h , z_h \cdot n_e \rangle_e & = \langle - \dgGrd^\ast z_h, y_h \rangle
	\end{alignat}
\end{subequations}
and
\begin{equation}
	\label{Eq:DiscdualV}	
\langle  \dgGrd y_h , z_h\rangle + \langle  \dgDiv z_h , y_h\rangle = - \sum_{e\in \mathcal E_h^\partial} \langle y_h , z_h \cdot n_e \rangle_e\,.
\end{equation}

Now, we are able to define a discrete counterpart $A_h: H_h \rightarrow H_h^\ast$ of the differential operator $A:D(A)\subset H \rightarrow H$ introduced in \eqref{Eq:HPS_31}. 

\begin{defi}[Discrete operator $A_h$]
\label{Def:Ah}
For the DG differential operators introduced in \eqref{Eq:dgOpHh} and \eqref{Eq:dgOpHhV}, respectively, the operator  $A_h: H_h \rightarrow H_h^\ast$ is defined by 
\begin{equation}
\label{Eq:DefAh}
A_h:= \begin{pmatrix}
	0 & -\dgDiv & 0 & 0 \\
	-\dgGrd & 0 & 0 & 0 \\ 
	0 & 0 & 0 & \dgdiv\\ 
	0 & 0 & \dggrd & 0
\end{pmatrix}\,,
\end{equation}
such that for $Y, Z\in H_h$ there holds that 
\begin{equation}
	\label{Eq:DefAh2}
	\langle A_h Y , Z \rangle   = - \langle \dgDiv Y_2, Z_{1}\rangle - \langle \dgGrd Y_1, Z_{2} \rangle +\langle \dgdiv Y_4, Z_{3} \rangle +  \langle \dggrd Y_3, Z_{4} \rangle \,.
\end{equation}	
\end{defi}

From \eqref{Eq:CondDualPtb3} and \eqref{Eq:DiscdualV} we conclude that for $Y\in H_h$ there holds that 
\begin{equation}
\label{Eq:DefAh3}
\langle A_h Y , Y \rangle   =  \sum_{e\in \mathcal E_h^\partial} \langle Y_1 , Y_2 \cdot n_e \rangle_e -  \sum_{e\in \mathcal E_h^\partial} \langle Y_3, Y_4 \cdot n_e \rangle_e\,.
\end{equation}
By \eqref{Eq:CondDualPtb0} and \eqref{Eq:CondDualPtb10}, the operator $A_h$ is not skew-selfadjoint on $H_h^{\operatorname{dg}}$, defined in \eqref{Eq:DefHhdg}, due to perturbations by boundary face integrals. Consequently, the inner product $\langle A_h Y , Y \rangle $ does no longer vanish as in the continuous case. However, the control of the latter term is essential for our analysis. Therefore, some correction term, defined in \eqref{Eq:FdPnc4} below, will be introduced in the fully discrete scheme. Finally, we note that skew-selfadjointness is preserved for the hybrid space $H_h^{\operatorname{hy}}$ of \eqref{Eq:DefHhhy}.


\subsection{Fully discrete problem of structure preserving nonconforming approximation}
\label{Subsec:NconfFD}

For the discretization of Problem~\ref{Prob:EP} in the space $\Hhnc$, defined in \eqref{Eq:DefHh} and \eqref{Eq:DefHhhydg}, we then consider the following family of fully discrete nonconforming approximation schemes. 

\begin{prob}[Fully discrete problem]
\label{Prob:FdPnc}
Let $\Hhnc$ be given by \eqref{Eq:DefHh} and \eqref{Eq:DefHhhydg}. For the operators $M_0$ and $M_1$ of \eqref{Eq:HPS_31}, $A_h$ of \eqref{Eq:DefAh}  and given data $F\in H^1_\nu(0,T;H)$ and $U_{0,h} \in \Hhnc$, where $U_{0,h} \in \Hhnc$ denotes an approximation of the initial value $U_0 \in H$ according to \eqref{Eq:HPS_30}, find $U_{\tau,h} \in Y_{\tau}^k(\Hhnc)$ with 
\begin{equation}
	\label{Eq:FdPnc1}
	\begin{aligned}
	Q_n[(\partial_t M_0 + M_1 + A_h)U_{\tau,h},V_{\tau,h}]_\nu & + Q_n[J_\partial(U_{\tau,h},V_{\tau,h})+J_\gamma(U_{\tau,h},V_{\tau,h})]_\nu\\[1ex]
	& + \langle M_0 \lsem U_{\tau,h} \rsem_{n-1}, V_{\tau,h}^{+,n-1}\rangle_H = Q_n[F,V_{\tau,h}]_\nu 
	\end{aligned}
\end{equation}
for all $V_{\tau,h} \in Y_{\tau}^k(\Hhnc)$ and $n\in \{1,2,\ldots,N\}$, where 
\begin{align}
\label{Eq:FdPnc4}
J_\partial(U_{\tau,h},V_{\tau,h}) & := - \sum_{e\in \mathcal E_h^\partial} \langle  U_{\tau,h}^{(2)} \cdot n_e, V_{\tau,h}^{(1)} \rangle_e +  \sum_{e\in \mathcal E_h^\partial} \langle U_{\tau,h}^{(4)}\cdot n_e, V_{\tau,h}^{(3)}  \rangle_e	\,, \\[1ex]
\label{Eq:FdPnc3}
J_\gamma(U_{\tau,h},V_{\tau,h}) & := \sum_{e\in \mathcal E_h^\partial}  \dfrac{1}{h} \Big(\gamma_1\langle U_{\tau,h}^{(1)}, V_{\tau,h}^{(1)}\rangle_e + \gamma_3 \langle U_{\tau,h}^{(3)}, V_{\tau,h}^{(3)} \rangle_e \Big)
\end{align}
and 
\begin{equation}
	\label{Eq:FdPnc2}
	\lsem U_{\tau,h} \rsem_{n-1} : = \left\{ \begin{array}{@{}ll} U_{\tau,h}^+(t_{n-1}) - U_{\tau,h}^-(t_{n-1})\,, & \text{for} \;\; n\in \{2,\ldots,N\}\,, \\[1ex]
		U_{\tau,h}^+(t_{n-1}) - U_{0,h}\,, & \text{for} \;\; n = 1\,,
	\end{array}\right.
\end{equation}
for $U_{\tau,h},V_{\tau,h}\in Y_{\tau}^k(\Hhnc)$ with $V_{\tau,h}^{+,n-1}:= V_{\tau,h}^+(t_{n-1})$; cf.~\eqref{Eq:DefLim}.
\end{prob}

\begin{rem}
\label{Rem:TMS}
\begin{itemize}
\item The algorithmic (or penalization) parameters $\gamma_i>0$, for $i\in \{1,2\}$, in \eqref{Eq:FdPnc3} have to be chosen sufficiently large; cf.~\cite{R08}. The  contribution $J_\gamma$, defined in \eqref{Eq:FdPnc3}, enforces the weak form of the homogeneous Dirichlet boundary conditions in \eqref{Eq:DefDA}. Moreover, in the error estimation given below it is essential for absorbing contributions from upper bounds of the error. 
	
\item In \eqref{Eq:FdPnc1}, the mathematical structure of the evolutionary problem \eqref{Eq:HPS_32} is essentially preserved, with the discrete operator $A_h$ replacing $A$. The perturbation of the skew-selfadjointness of $A_h$,  resulting from \eqref{Eq:CondDualPtb0} and \eqref{Eq:CondDualPtb10}, is captured in the analysis below by the additional (boundary) correction $Q_n[J_\partial(U_{\tau,h},V_{\tau,h})]_\nu$ along with the penalization induced by $Q_n[J_\gamma(U_{\tau,h},V_{\tau,h})]_\nu$ . 
	
\item Problem \ref{Prob:FdPnc} yields a global in time formulation. For computations of space-time finite element discretizations we propose using a temporal test basis that is supported on the subintervals $I_n$; cf.~\cite{AB23,ABMS23}. Then, a time marching process is obtained. For Problem~\ref{Prob:FdPnc}, this amounts to assuming that the trajectory $U_{ \tau,h}$ has been computed before for all $t\in [0,t_{n-1}]$, starting with an approximation $U_{\tau,h}(t_0) := U_{0,h}$ of $U_0\in D(A)$. On $I_n=(t_{n-1},t_n]$, for given $U_{\tau,h}(t_{n-1})\in \Hhnc$ we consider then finding $U_{\tau,h}\in \mathbb P_k(I_n,\Hhnc)$ such that \eqref{Eq:FdPnc1} is satisfied for all $V_{\tau,h}\in \mathbb P_k(I_n,\Hhnc)$.

\item In \eqref{Eq:FdPnc1}, there holds that  
\begin{equation}
\label{Eq:DivFreeFEM}
\langle A_h U_{\tau,h}, V_{\tau,h}\rangle_H = \Bigg\langle \begin{pmatrix} -\dgDiv U^{(2)}_{\tau,h}\\[0.5ex] -\dgGrd U^{(1)}_{\tau,h}\\[0.5ex] \dgdiv U^{(4)}_{\tau,h}\\[0.5ex] \dggrd U^{(3)}_{\tau,h} \end{pmatrix}, \begin{pmatrix} V^{(1)}_{\tau,h} \\[0.5ex] V^{(2)}_{\tau,h} \\[0.5ex] V^{(3)}_{\tau,h} \\[0.5ex] V^{(4)}_{\tau,h} \end{pmatrix}\Bigg\rangle_H =  \Bigg\langle \begin{pmatrix} U^{(2)}_{\tau,h}\\[0.5ex] -\dgGrd U^{(1)}_{\tau,h}\\[0.5ex] U^{(4)}_{\tau,h}\\[0.5ex] \dggrd U^{(3)}_{\tau,h} \end{pmatrix}, \begin{pmatrix} -\dgDiv^\ast V^{(1)}_{\tau,h} \\[0.5ex] V^{(2)}_{\tau,h} \\[0.5ex] \dgdiv^\ast V^{(3)}_{\tau,h} \\[0.5ex] V^{(4)}_{\tau,h} \end{pmatrix}\Bigg\rangle_H\,.
\end{equation}
By \eqref{Eq:CondDualPtb1} and \eqref{Eq:CondDualPtb11}, the operators $-\dgDiv^\ast$ and $\dgdiv^\ast$ in \eqref{Eq:DivFreeFEM} are transformed into the DG gradients $\dgGrd$ and $\dggrd$, respectively, applied to the test functions, and additional sums of boundary face integrals. This can be exploited in the assembly process and error analysis.
\end{itemize}
\end{rem}

\begin{thm}[Well-posedness of fully discrete problem]
\label{Thm:WPFdPd}
There exists a unique solution $U_{\tau,h} \in Y_{\tau,h}^k(\Hhnc)$ of Problem \ref{Prob:FdPnc}. 
\end{thm}

\begin{mproof}
The proof follows the ideas of \cite[Proof of Prop.\ 3.2]{FTW19}. To keep this work self-contained, and due to adaptations of the proof required by the perturbation of the skew-selfadjointness, we present it briefly. Since Problem \ref{Prob:FdPnc} is finite dimensional, it suffices to prove uniqueness of solutions to \eqref{Eq:FdPnc1} for $n\in \{1,\ldots,N\}$. The existence of solutions then directly follows from their uniqueness. By means of the first of the items in Rem.~\ref{Rem:TMS} and an induction argument, it suffices to prove the uniqueness of solutions to \eqref{Eq:FdPnc1} on a fixed subinterval $I_n$. For this, let $\widetilde U_{\tau,h}\in \mathbb P_k(I_n,\Hhnc)$ and $\widehat U_{\tau,h}\in \mathbb P_k(I_n,\Hhnc)$ be two solutions of \eqref{Eq:FdPnc1}. Then, their difference $ U_{\tau,h}:=\widetilde U_{\tau,h}-\widehat U_{\tau,h}$ satisfies for all $V_{\tau,h}\in \mathbb P_k(I_n,\Hhnc)$ that
\begin{equation}
	\label{Eq:Uniq01}
	\begin{aligned}
	Q_n[(\partial_t M_0 + M_1 + A_h)U_{\tau,h},V_{\tau,h}]_\nu & + Q_n[J_\partial(U_{\tau,h},V_{\tau,h})+ J_\gamma(U_{\tau,h},V_{\tau,h})]_\nu\\[1ex]
	 & + \langle M_0 U_{\tau,h}^+(t_{n-1})\,. V_{\tau,h}^{+,n-1}\rangle_H = 0\,.
	\end{aligned}
\end{equation}

Next, we recall an argument of \cite[Proof of Prop.\ 3.2]{FTW19}. We note that 
\begin{equation*}
	\partial_t: \mathbb P_k(I_n,H) \rightarrow \mathbb  P_k(I_n,H)\,, \;\; w_\tau \mapsto \partial_t w_\tau 
\end{equation*}
and 
\begin{equation*}
	\delta_{n-1}: \mathbb  P_k(I_n,H) \rightarrow H \,, \;\; w_\tau \mapsto \delta_{n-1} w_\tau := w^+_\tau(t_{n-1})
\end{equation*}
are bounded linear operators with respect to the norm of $\mathbb P_k(I_n,H)$ induced by the inner product \eqref{Def:IPPk}. Consequently, the mapping 
\begin{equation*}
\mathbb P_k(I_n,H) \rightarrow \R \,, \;\; w_\tau \mapsto \langle z , \delta_{n-1} w_\tau\rangle_H
\end{equation*}
is linear and bounded for each $z \in H$. Then, by the Riesz representation theorem there exists a unique $\Psi_\tau(z)\in \mathbb P_k(I_n;H)$ such that 
\begin{equation}
\label{Eq:Uniq02}
	\langle \Psi_\tau(z),w_\tau\rangle_{\nu,n} = \langle z, \delta_{n-1} w_\tau\rangle_H\,. 
\end{equation}
The mapping $\Psi_\tau : H \rightarrow \mathbb P_k(I_n,H)$ is linear and bounded, since for $z\in H$ there holds that 
\begin{equation*}
\|\Psi_\tau(z)\|^2_{\nu,n} = 	\langle \Psi_\tau(z),\Psi_\tau(z)\rangle_{\nu,n} = \langle z , \delta_{n-1} \Psi_\tau(z)\rangle_H \leq \| z \|_H \|\delta_{n-1}\| \| \Psi_\tau(z) \|_{\nu,n} \,.
\end{equation*}
Now, using integration by parts along with \eqref{Eq:Uniq02}, we have or all $v_\tau\in \mathbb P_k(I_n;H)$ that 
\begin{align}
\nonumber
\langle \partial_t M_0 v_\tau, v_\tau \rangle_{\nu,n} & = \dfrac{1}{2} \langle \partial_t M_0 v_\tau, v_\tau \rangle_{\nu,n} + \frac{1}{2} \int_{t_{n-1}}^{t_n} \langle M_0 \partial_t v_\tau(t),  v_\tau(t)\rangle_H \operatorname{e}^{-2\nu (t-t_{n-1})} \ud t\\[1ex]
\nonumber
& = \dfrac{1}{2} \langle \partial_t M_0 v_\tau, v_\tau \rangle_{\nu,n} - \dfrac{1}{2} \int_{t_{n-1}}^{t_n} \langle  M_0 \partial_t v_\tau(t),  v_\tau(t)\rangle_H \operatorname{e}^{-2\nu (t-t_{n-1})} \ud t\\[1ex]
\nonumber
& \quad + \nu \int_{t_{n-1}}^{t_n} \langle  M_0  v_\tau(t), v_\tau(t)\rangle_H \operatorname{e}^{-2\nu (t-t_{n-1})} \ud t\\[0.75ex]
\nonumber
& \quad  + \frac{1}{2} \langle  M_0  v_\tau(t_n),v_\tau(t_n)\rangle_H \operatorname{e}^{-2\nu \tau_n} - \frac{1}{2} \langle v^+_\tau(t_{n-1}), M_0 v^+_\tau(t_{n-1})\rangle_H\\[1ex]
\label{Eq:Uniq04}
& \geq \nu \langle M_0 v_\tau,v_\tau\rangle_{\nu,n} - \frac{1}{2} \langle \Psi_\tau(M_0 \delta_{n-1} v_\tau), v_\tau\rangle_{\nu,n}\,.
\end{align}

Using \eqref{Eq:Uniq02}, we rewrite \eqref{Eq:Uniq01} as  
\begin{equation}
\label{Eq:Uniq03}
\begin{aligned}
Q_n[(\partial_t M_0 + M_1 + A_h)U_{\tau,h},V_{\tau,h}]_\nu & + Q_n[J_\partial(U_{\tau,h},V_{\tau,h})+J_\gamma(U_{\tau,h},V_{\tau,h})]_\nu\\[1ex] 
& + \langle \Psi_\tau(M_0 \delta_{n-1} U_{\tau,h}), V_{\tau,h} \rangle_H = 0\,.
\end{aligned}
\end{equation}
In \eqref{Eq:Uniq03}, we choose $V_{\tau,h}=U_{\tau,h}$. By \eqref{Eq:DefAh3} along with \eqref{Eq:FdPnc4} we have for $Z \in Y_h$ that  
\begin{align}
	\nonumber 
	\langle A_h Z, Z\rangle + J_\partial(Z,Z) & =   \sum_{e\in \mathcal E_h^\partial} \langle Z_1 , Z_2 \cdot n_e \rangle_e -  \sum_{e\in \mathcal E_h^\partial} \langle Z_3, Z_4 \cdot n_e \rangle_e\\[1ex]
	\label{Eq:Uniq05}
	& \quad - \sum_{e\in \mathcal E_h^\partial} \langle Z_2 \cdot n_e, Z_1 \rangle_e +   \sum_{e\in \mathcal E_h^\partial} \langle Z_4 \cdot n_e , Z_3\rangle_e  = 0\,.
\end{align}
Now, from \eqref{Eq:Uniq03} we deduce by \eqref{Eq:Uniq05} and the nonnegativity of $J_\gamma(U_{\tau,h},U_{\tau,h})$ given by \eqref{Eq:FdPnc3} that 
\begin{align}
\nonumber
0 & = Q_n[(\partial_t M_0 + M_1 + A_h)U_{\tau,h},U_{\tau,h}]_\nu + Q_n[J_\partial(U_{\tau,h},U_{\tau,h})+J_\gamma(U_{\tau,h},U_{\tau,h})]_\nu\\
& \quad + \langle \Psi_\tau(M_0 \delta_{n-1} U_{\tau,h}), U_{\tau,h} \rangle_H \\[1ex]
\nonumber
& \geq  \langle \partial_t M_0 U_{\tau,h}, U_{\tau,h} \rangle_{\nu,n} + \langle M_1 U_{\tau,h}, U_{\tau,h} \rangle_{\nu,n} + \langle \Psi_\tau(M_0 \delta_{n-1} U_{\tau,h}), U_{\tau,h}\rangle_{\nu,n}\\[1ex]
\nonumber
& \geq  \langle (\nu M_0 + M_1)U_{\tau,h} , U_{\tau,h}\rangle_{\nu,n}  +  \frac{1}{2} \langle \Psi_\tau(M_0 \delta_{n-1} U_{\tau,h}), U_{\tau,h} \rangle_{\nu,n}\\[1ex]
\label{Eq:Uniq06}
& \geq \gamma \langle U_{\tau,h}, U_{\tau,h} \rangle_{\nu,n}\,,
\end{align}
where the nonnegativity of $\langle \Psi_\tau(M_0 \delta_{n-1} U_{\tau,h}), U_{\tau,h} \rangle_{\nu,n}$ is ensured by 
\begin{equation*}
\langle \Psi_\tau(M_0 \delta_{n-1} U_{\tau,h}), U_{\tau,h} \rangle_{\nu,n} = \langle M_0  U_{\tau,h}(t_{n-1}^+),U_{\tau,h}(t_{n-1}^+) \rangle_H \geq 0\,.
\end{equation*}	
The latter inequality follows from the assumption \eqref{Eq:WP0}. From \eqref{Eq:Uniq06} we directly conclude the uniqueness of solutions to \eqref{Eq:FdPnc1} and, thereby, the assertion of this lemma.
\end{mproof}

\section{Error estimation for structure preserving nonconforming approximation}
\label{Subsec:ErrEst}

Here we prove an error estimate for the solution $U_{\tau,h}$ of Problem~\ref{Prob:FdPnc}. For brevity, the proof is done only for the full DG approximation in space, corresponding to the choice $\Hhnc=H_h^{\text{dg}}$ in \eqref{Eq:DefHh} with $r \in \N_0$ in \eqref{Eq:DefHhdg}. The adaptation of the proof to the hybrid case $\Hhnc=H_h^{\text{hy}}$ in \eqref{Eq:DefHh} is straightforward.  

\begin{thm}[Error estimate for the fully discrete problem]
	\label{Thm:ErrEstSTnc}
	Let $H$ be defined by \eqref{Eq:DefH}. For the solution $U$ of Problem~\ref{Prob:EP} suppose that the regularity condition
	\begin{equation}
		\label{Eq:Reg0}
		U  \in H^{k+3}_\nu(\R;H)\cap H^2_\nu\big(\R;H^{r+1}(\Omega)^{(d+1)^2} \big)
	\end{equation}
	is satisfied. Let the discrete initial $U_{0,h}\in \Hhnc$ in Problem~\ref{Prob:FdPnc} be chosen such that $\|U_0-U_{0,h}\|\leq c h^{2r}$ holds.  Then, for the numerical solution $U_{\tau,h}$ of Problem \ref{Prob:FdPnc} we have the error estimate that
	\begin{equation}
		\label{Eq:ErrEstSTnc01}
		\sup_{t\in [0,T]} \langle M_0 (U-U_{\tau,h})(t), (U-U_{\tau,h})(t) \rangle+  \operatorname{e}^{2\nu T} \| U- U_{\tau,h}\|_{\tau,\nu}^2 \leq C (1+T) \operatorname{e}^{2\nu T} (\tau^{2(k+1)} + h^{2r})\,.
	\end{equation}
\end{thm}

\begin{mproof}
	We split the error $U - U_{\tau,h}$ into the two parts 
	\begin{equation}
		\label{Eq:ErrEstSTnc02}
		U - U_{\tau,h} = Z + E_{\tau,h} \qquad \text{with}\quad Z:= U - I_\tau \Pi_h  U\,,  \quad  E_{\tau,h}:= I_\tau \Pi_h  U - U_{\tau,h} \,,
	\end{equation}
	where $I_\tau$ and $\Pi_h $ are defined in \eqref{Def:LIO} and \eqref{Def:Pih}, respectively. The errors $Z$ and $E_{\tau,h}$ are estimated  in Lem.~\ref{Lem:ErrProj} to \ref{Lem:ErrOrderEst} below. By means of the triangle inequality, the splitting \eqref{Eq:ErrEstSTnc02} along with these lemmas then proves  \eqref{Eq:ErrEstSTnc01}. 
\end{mproof}

For the error $Z$ in \eqref{Eq:ErrEstSTnc02} there holds the following estimate.  
\begin{lem}[Estimation of error $Z$]
	\label{Lem:ErrProj}
	For $s\in \{0,\ldots, r+1\}$ let  $U\in H_\nu^{k+2}(\R;H)\cap H_\nu^1(\R;H^{s}(\Omega)^{(d+1)^2})$ be satisfied. For the error $Z=U- I_\tau \Pi_h  U$ there holds that 
	\begin{equation}
		\label{Eq:ErrProj0}
		\sup_{t\in[0,T]} \langle M_0 Z(t), Z(t) \rangle_H +  \operatorname{e}^{2\nu T} \| Z \|_{\tau,\nu}^2 \leq C (1+T) \operatorname{e}^{2\nu T}  (\tau^{2(k+1)} + h^{2s})\,.
	\end{equation}
\end{lem}

\begin{mproof}
	We split the error $Z$ into the two parts 
	\begin{equation}
		\label{Eq:ErrProj1}
		Z = U -  I_\tau \Pi_h U  = (U - \Pi_h U) +  (\Pi_h U -  I_\tau \Pi_h U) \,.
	\end{equation}	
	From \eqref{Eq:ErrProj1} along with the commutativity of $\Pi_h$ and $I_\tau$ and the boundedness of $M_0$ and $\Pi_h$ we get that 
	\begin{equation}
		\label{Eq:ErrProj2}
		\begin{aligned}
			\langle M_0 Z(t), Z(t) \rangle_H	& \leq C \|  Z(t) \|_H^2 = C \|  U(t) -  I_\tau \Pi_h U(t) \|_H^2\\[1ex]
			& \leq C \Big(\|  U(t) -  \Pi_h U(t) \|_H^2 + \|  U(t) -  I_\tau U(t) \|_H^2\Big)\,.
		\end{aligned}
	\end{equation}
	Using \eqref{Eq:ErrLI} and \eqref{Eq:IOPih}, we obtain from \eqref{Eq:ErrProj2} that 
	\begin{equation}
		\label{Eq:ErrProj3}
		\begin{aligned}
			\langle M_0 Z(t), Z(t) \rangle_H	& \leq C \Big(\tau^{2(k+1)} \sup_{t\in I_n} \| \partial_t^{k+1} U \|_H^2+ h^{2s}  \sup_{t\in I_n} \| U \|_{H^s(\Omega)}^2 \Big)\,, \quad \text{for} \;\; t\in I_n\,.
		\end{aligned}
	\end{equation}
	We note that the norms on the right-hand side of \eqref{Eq:ErrProj3} remain finite under the assumptions made about $U$; cf.\ \eqref{Eq:SET}. This shows the first of the  estimates in \eqref{Eq:ErrProj0}. By the definition of $\| \cdot \|_{\tau,\nu}$ in \eqref{Eq:TdN}, the second of the estimates in \eqref{Eq:ErrProj0} follows from \eqref{Eq:ErrProj2} along with \eqref{Eq:ErrLI} and \eqref{Eq:IOPih}. 
\end{mproof}

For the error $E_{\tau,h}$ in \eqref{Eq:ErrEstSTnc02} there holds the following estimate.  

\begin{lem}[Estimation of error $E_{\tau,h}$]
	\label{Lem:ErrDisc}
	For the error $E_{\tau,h}= I_\tau \Pi_h  U - U_{\tau,h}$ there holds that 
	\begin{equation}
		\label{Eq:ErrDisc0}
		\begin{aligned}
			\langle M_0 E_{\tau,h}^-(t_{N})  , E_{\tau,h}^-(t_{N})\rangle_H & +  \operatorname{e}^{2\nu T} (\| E_{\tau,h} \|_{\tau,\nu}^2  +  | J_\gamma (E_{\tau,h},E_{\tau,h})| _{\tau,\nu}^2)\\[1ex] 
			& \leq C \operatorname{e}^{2\nu T} \Big(\langle M_0 E_{\tau,h}^-(t_0), E_{\tau,h}^-(t_0)\rangle_H+ \| \partial_t M_0 (U - \hItau  \Pi_h U)\|_{\tau,\nu}^2 \\[1ex]
			& \quad + \| M_1 Z \|_{\tau,\nu}^2  + \| A_h Z \|_{\tau,\nu}^2+   | J_\gamma (Z,Z)|_{\tau,\nu}^2 + |J_\partial^{n} (Z,Z)|_{\tau,\nu}^2\\[1ex]
			& \quad  + T\max_{1\leq n\leq N} \Big\{ \|M_0(\Pi_h U^+(t_{n-1}) - I_\tau \Pi_h U^+(t_{n-1})) \|_H^2 \operatorname{e}^{-2\nu t_{n-1}} \Big\}\Big)\,,
		\end{aligned}	
	\end{equation}	
	where 
	\begin{equation}
		\label{Eq:ErrDisc00}
	 J_\partial^{n}(Z,Z):=\sum_{e\in \mathcal E_h^\partial} h \, \big(\langle Z_2\cdot n,Z_2\cdot n\rangle_e +  \langle Z_4\cdot n,Z_4\cdot n\rangle_e\big)\,.
	 \end{equation}
\end{lem}

\begin{mproof}
	Essentailly, the proof follows \cite[Thm.\ 3.8]{FTW19} for the semidiscretization in time. In \cite[Thm.\ 3.8]{FTW19} the skew-selfadjointness of the continuous operator is a key ingredient for proving the error estimate. To keep this work self-contained, we summarize in the appendix the proof of \eqref{Eq:ErrDisc0} for the setting of Problem~\ref{Lem:ErrDisc} and the perturbed skew-selfadjointness of $A_h$,  depicted by  Lem.~\ref{Lem:CondDualIBV} and \eqref{Eq:CondDualPtb10}.  
\end{mproof}

For the error $E_{\tau,h}$ in \eqref{Eq:ErrEstSTnc02} we further have the following improved estimate. 

\begin{lem}[Improved estimation of error $E_{\tau,h}$]
	\label{Lem:ErrDiscInf}
	For the error $E_{\tau,h}$ in \eqref{Eq:ErrEstSTnc02} there holds that 
	\begin{equation}
		\label{Eq:ErrDiscInf0}
		\begin{aligned}
			\sup_{t\in [0,T]} \langle M_0 E_{\tau,h}(t) & , E_{\tau,h}(t) \rangle_H \leq C  \Big(\langle M_0 E_{\tau,h}^-(t_0), E_{\tau,h}^-(t_0)\rangle_H + \| \partial_t M_0 (U - \hItau  \Pi_h U)\|_{\tau,\nu}^2 \\[1ex] 
			& + \| M_1 Z \|_{\tau,\nu}^2+ \| A_h Z \|_{\tau,\nu}^2 +   | J_\gamma (Z,Z)|_{\tau,\nu}^2 + |J_\partial^{n} (Z,Z)|_{\tau,\nu}^2\\[1ex]
			&  + T\max_{1\leq n\leq N} \Big\{ \|M_0(\Pi_h U^+(t_{n-1}) - I_\tau \Pi_h U^+(t_{n-1})) \|_H^2 \operatorname{e}^{-2\nu t_{n-1}} \Big\}\Big)\,.
		\end{aligned}	
	\end{equation}	
\end{lem}

\begin{mproof}
	The proof follows the lines of \cite[Thm.\ 3.12]{FTW19} for the semidiscretization in time. Again, to keep this work self-contained, we summarize the proof for the setting of Problem~\ref{Prob:FdPnc} in the appendix.  
\end{mproof}

Next, we estimate the terms on the right-hand side of \eqref{Eq:ErrDisc0} and \eqref{Eq:ErrDiscInf0}, respectively, one by one. 

\begin{lem}
	\label{Lem:ErrOrderEst}
	For $s\in \{1,\ldots, r+1\}$ let  $U\in H_\nu^{k+3}(\R;H)\cap H_\nu^2(\R;H^{s}(\Omega)^{(d+1)^2})$ be satisfied. With $Z=U- I_\tau \Pi_h  U$ there holds that 
	\begin{subequations}
		\label{Eq:ErrOrderEst0}
		\begin{alignat}{2}
			\label{Eq:ErrOrderEst1}
			\| \partial_t M_0 (U - \hItau  \Pi_h U)\|_{\tau,\nu}^2 +  \| M_1 Z \|_{\tau,\nu}^2&  \leq C T  (\tau^{2(k+1)}+h^{2s})\,,\\[1ex]
			\label{Eq:ErrOrderEst2}
			\| A_h Z \|_{\tau,\nu}^2 +   | J_\gamma (Z,Z)|_{\tau,\nu}^2 + |J_\partial^{n} (Z,Z)|_{\tau,\nu}^2 & \leq C T h^{2(s-1)}\,,\\[1ex]
			\label{Eq:ErrOrderEst3}
			\max_{1\leq n\leq N} \Big\{ \|M_0(\Pi_h U^+(t_{n-1}) - I_\tau \Pi_h U^+(t_{n-1})) \|_H^2 \operatorname{e}^{-2\nu t_{n-1}} \Big\}\Big) & \leq C \tau^{2(k+1)}\,.	
		\end{alignat}
	\end{subequations}
\end{lem}

\begin{mproof}
	Using a splitting as in \eqref{Eq:ErrProj1} and commutation properties of the operators shows that   
	\begin{equation}
		\label{Eq:ErrOrderEst4}
		\begin{aligned}
			\| \partial_t M_0 (U(t) - \hItau  \Pi_h U(t))\|_H  & \leq \| M_0 \Pi_h \partial_t (U(t) - \hItau  U(t) )\|_H\\[1ex] 
			& \quad +  \| M_0 (\partial_t U(t) - \Pi_h \partial_t U(t)) \|_H\,.
		\end{aligned}
	\end{equation}	
	By \eqref{Eq:ErrLI0c} and \eqref{Eq:IOPih} along with the boundedness of $M_0$ and $\Pi_h$ we conclude from \eqref{Eq:ErrOrderEst4} that 
	\begin{equation}
		\label{EE:ErrOrderEst5}
		\| \partial_t M_0 (U(t) - \hItau  \Pi_h U(t))\|_H  \leq C \Big(\tau^{k+1} \sup_{t\in I_n} \| \partial_t^{k+2} U \|_H+ h^{s}  \sup_{t\in I_n} \| \partial_t U \|_{H^{s}(\Omega)} \Big)\,.
	\end{equation}	
	Recalling the definition of the time-discrete norm $\|\cdot \|_{\tau,\nu}$ in \eqref{Eq:TdN}, this directly proves the first of the bounds in \eqref{Eq:ErrOrderEst1}. The estimate of $ \| M_1 Z \|_{\tau,\nu}$ in \eqref{Eq:ErrOrderEst1} and inequality \eqref{Eq:ErrOrderEst3} follow similarly. 
	
	It remains to prove \eqref{Eq:ErrOrderEst2}. By \eqref{Def:LIOc}, the interpolation operator $\hItau$ acts as the identity in the Gauss--Radau points $t_{n,\mu}$. Thus, for $\mu=0,\ldots,k$ we get that 
	\begin{subequations}
		\label{EE:ErrOrderEst6}
		\begin{alignat}{2}
		\label{EE:ErrOrderEst6a}
		\| A_h Z(t_{n,\mu}) \|_H  & = \| A_h (U -\Pi_h U)(t_{n,\mu})\|_H\,, \\[1ex]
		\label{EE:ErrOrderEst6b}
		 | J_\gamma (Z(t_{n,\mu}),Z(t_{n,\mu}))|_{\tau,\nu} & =  | J_\gamma ((U -\Pi_h U)(t_{n,\mu}),(U -\Pi_h U)(t_{n,\mu}))|_{\tau,\nu}  \,, \\[1ex] 
		\label{EE:ErrOrderEst6c}
		 |J_\partial^{n} (Z,Z)|_{\tau,\nu} & = |J_\partial^{n} ((U -\Pi_h U)(t_{n,\mu}),(U -\Pi_h U)(t_{n,\mu}))|_{\tau,\nu}\,.
		\end{alignat}
	\end{subequations}	
	Putting $\Theta^n_\mu := (U - \Pi_h U)(t_{n,\mu})$, using \eqref{Eq:DivFreeFEM} and recalling the duality relations \eqref{Eq:CondDualPtb1} and \eqref{Eq:CondDualPtb11}, we get for the right-hand side of \eqref{EE:ErrOrderEst6a} that 
	\begin{equation}
		\label{EE:ErrOrderEst7}
		\begin{aligned}
			\langle A_h  \Theta^n_\mu, Y_h \rangle  & = \langle \Theta^n_{\mu,2}, -\dgDiv^\ast Y_{h,1}\rangle - \langle \dgGrd \Theta^n_{\mu,1}, Y_{h,2}\rangle \\[1ex]
			& \quad + \langle \Theta^n_{\mu,4}, \dgdiv^\ast Y_{h,3}\rangle + \langle \dggrd \Theta^n_{\mu,3}, Y_{h,4} \rangle\\[1ex]
			& =  \langle \Theta^n_{\mu,2}, \dgGrd Y_{h,1}\rangle  + \sum_{e\in \mathcal E_h^\partial} \langle \Theta^n_{\mu,2} \cdot n, Y_{h,1} \rangle_e - \langle \dgGrd \Theta^n_{\mu,1}, Y_{h,2}\rangle \\[1ex]
			& \quad - \langle \Theta^n_{\mu,4}, \dggrd Y_{h,3}\rangle  - \sum_{e\in \mathcal E_h^\partial} \langle \Theta^n_{\mu,4}\cdot n, Y_{h,3} \rangle_e + \langle \dggrd \Theta^n_{\mu,3}, Y_{h,4} \rangle
		\end{aligned}
	\end{equation}	
	for all $Y_h \in H_h^{\text{dg}}$. Next, we bound the right-hand side of \eqref{EE:ErrOrderEst7} term by term. We start with the last term in \eqref{EE:ErrOrderEst7}. By \eqref{Eq:dggrdHh}, defining $\dggrd$, and the Cauchy--Schwarz inequality it follows that 
	\begin{align}
		\nonumber
		\langle \dggrd  \Theta^n_{\mu,3} , Y_{h,4} \rangle  & = \langle \operatorname{grad}_h \Theta^n_{\mu,3}, Y_{h,4}\rangle - \sum_{e\in \mathcal E_h^i} \langle \lsem \Theta^n_{\mu,3}\rsem, \ldblbrace  Y_{h,4} \rdblbrace \cdot n_e \rangle_e  - \sum_{e\in \mathcal E_h^\partial } \langle  \Theta^n_{\mu,3},  Y_{h,4} \cdot n_e \rangle_e  \\[1ex]
		\label{EE:ErrOrderEst8}
		& \leq \sum_{K\in \mathcal T_h} \|\nabla \Theta^n_{\mu,3}\|_{L^2(K)} \| Y_{h,4}\|_{L^2(K)} + \sum_{e\in \mathcal E_h^i} \| \lsem \Theta^n_{\mu,3} \rsem \|_{L^2(e)} \|  \ldblbrace Y_{h,4} \rdblbrace \cdot n_e\|_{L^2(e)} \\[1ex]
		\nonumber
		& \quad + \sum_{e\in \mathcal E_h^\partial} \| \Theta^n_{\mu,3} \|_{L^2(e)} \|  Y_{h,4} \cdot n_e \|_{L^2(e)}\,.
	\end{align}	
	Using \eqref{Eq:IOPih} with $m=1$, \eqref{Eq:IOPih21} and the inverse relation (cf.~\cite[Lem.~1.46]{PE12}) 
	\begin{equation}
		\label{EE:ErrOrderEst9}
		h_K^{1/2} \| w_h \|_{L^2(e)} \leq C_{\operatorname{inv}} \| w_h \|_{L^2(K)}\,, \quad \text{for } e\subset K\,, \;\; w\in \mathbb Q_r^d\,, 
	\end{equation}	
	we obtain from \eqref{EE:ErrOrderEst8} that 
	\begin{equation}
		\label{EE:ErrOrderEst10}
		\langle \dggrd \Theta^n_{\mu,3}, Y_{h,4} \rangle \leq C h^{s-1} \| U \|_{H^{s}(\Omega)} \Big(\sum_{K\in \mathcal T_h} \|Y_{h,4} \|_{L^2(K)}^2\Big)^{1/2}\,.
	\end{equation}	
	For the fourth and fifth term on right-hand side of \eqref{EE:ErrOrderEst7} we have that 
	\begin{equation*}
		\begin{aligned}
		    \langle \Theta^n_{\mu,4}, \dggrd Y_{h,3}\rangle  +  \sum_{e\in \mathcal E_h^\partial} \langle \Theta^n_{\mu,4}\cdot n, Y_{h,3} \rangle_e= \langle \Theta^n_{\mu,4}, \operatorname{grad}_h  Y_{h,3}\rangle - \sum_{e\in \mathcal E_h^i} \langle \ldblbrace  \Theta^n_{\mu,4} \rdblbrace \cdot n_e, \lsem Y_{h,3}\rsem \rangle_e  \,.
		\end{aligned}
	\end{equation*}
From this, we find that 
	\begin{equation}
		\label{EE:ErrOrderEst11}	
			\langle \Theta^n_{\mu,4}, \dggrd Y_{h,3}\rangle 	 \leq 	\sum_{K\in \mathcal T_h} \| \Theta^n_{\mu,4}\|_{L^2(K)} \| \nabla Y_{h,3}\|_{L^2(K)} + C \sum_{e\in \mathcal E_h^i} \| \ldblbrace  \Theta^n_{\mu,4} \rdblbrace \cdot n_e\|_{L^2(e)} \|\lsem Y_{h,3}\rsem \|_{L^2(e)} \,. 
	\end{equation}
	Using \eqref{Eq:IOPih} with $m=0$, bounding $\| \nabla Y_{h,3}\|_{L^2(K)}$ by the the $H^1$--$L^2$ inverse inequality and applying \eqref{Eq:IOPih21} and \eqref{EE:ErrOrderEst9}, we deduce from \eqref{EE:ErrOrderEst11} that 
	\begin{equation}
		\label{EE:ErrOrderEst12}
		\langle \Theta^n_{\mu,4}, \dggrd Y_{h,3}\rangle  \leq C h^{s-1} \| U \|_{H^{s}(\Omega)} \Big(\sum_{K\in \mathcal T_h} \|Y_{h,3} \|_{L^2(K)}^2\Big)^{1/2}\,.
	\end{equation}
	The first three terms on right-hand side of \eqref{EE:ErrOrderEst7} can be treated similarly. Since 
	\begin{equation*}
		\| A_h W \|_H = \sup_{Y_h \in H_h^{\text{dg}}\backslash\{0\}} \dfrac{\langle A_h W, Y_h \rangle_H}{\| Y_h \|_H} \,, \quad \text{for } \; W\in D(A)+H_h^{\text{dg}}\,,
	\end{equation*}
	combining \eqref{EE:ErrOrderEst7} with \eqref{EE:ErrOrderEst10} and \eqref{EE:ErrOrderEst12} and their counterparts for the first and the second term on right-hand side of \eqref{EE:ErrOrderEst7} proves for \eqref{EE:ErrOrderEst6} that 
	\begin{equation}
		\label{EE:ErrOrderEst13}
		\| A_h (U(t_{n,\mu}) -\Pi_h U(t_{n,\mu}))\|_H^2 \leq C h^{2(s-1)}\,, \quad \text{for }\; \mu =1,\ldots,k \,.	
	\end{equation} 
	 Applying the temporal quadrature formula \eqref{Eq:GRF} to \eqref{EE:ErrOrderEst13}, summing up the resulting inequality from $n=1$ to $N$ and recalling the definition in \eqref{Eq:TdN} yields the error bound for $\|A_h Z\|_{\tau,\nu}^2$ in \eqref{Eq:ErrOrderEst2}. The bound for $| J_\gamma (Z,Z)|_{\tau,\nu}^2$ follows similarly. Using \eqref{Eq:IOPih21} we get that 
	 \begin{equation}
		\label{EE:ErrOrderEst14}
	 	J_\gamma(\Theta_\mu^n,\Theta_\mu^n) \leq \dfrac{C}{h} h^{2s-1} \big(\| U^{(1)} \|_{H^{s}(\Omega)}^2 + \| U^{(3)} \|_{H^{s}(\Omega)}^2\big) \,.
	\end{equation}	
	 Applying the temporal quadrature formula \eqref{Eq:GRF} to \eqref{EE:ErrOrderEst14}, summing up the resulting inequality from $n=1$ to $N$ and recalling the definition in \eqref{Eq:TdAV} yields the error bound for $| J_\gamma (Z,Z)|_{\tau,\nu}^2$ in \eqref{Eq:ErrOrderEst2}. The term $|J_\partial^{n} (Z,Z)|_{\tau,\nu}^2$ in \eqref{Eq:ErrDisc00} is bounded by the same arguments. This completes the proof of \eqref{Eq:ErrOrderEst0}. 
\end{mproof}

Thm.\ \ref{Thm:ErrEstSTnc} proves convergence of the error measured in the time-mesh dependent norm $\|\cdot \|_{\tau,\nu}$ defined in \eqref{Eq:TdN}. Using a result of \cite[Thm.\ 2.5]{F23}, convergence with respect to the norm $\|\cdot \|_{\nu}$ induced by the inner product \eqref{Def:SPHnu} of the continuous function space $H_\nu(\R;H)$ can still be ensured. 

\begin{cor}
	\label{Lem:L2L2Err}
	Under the assumptions of Thm.\ \ref{Thm:ErrEstSTnc}, there holds that 	
	\begin{equation}
		\label{Eq:L2L2Err0}
		\| U- U_{\tau,h}\|_{\nu}^2 \leq C  (1+T) \operatorname{e}^{2\nu T} (\tau^{2(k+1)} + h^{2r})\,.
	\end{equation}
\end{cor} 

\begin{mproof}
	We split the error into the parts 
	\begin{equation}
		\label{Eq:L2L2Err1}
		\| U- U_{\tau,h}\|_{\nu} \leq  \| U- I_\tau U\|_{\nu} + \| I_\tau U - U_{\tau,h}\|_{\nu} \,.
	\end{equation}
	For the first of the terms on the right-hand side of \eqref{Eq:L2L2Err1}, we deduce from the interpolation error estimate proved in \cite[Thm.\ 2.5]{F23} that 
	\begin{equation}
		\label{Eq:L2L2Err2}
		\| U- I_\tau U\|_{\nu}^2 \leq  C \tau^{2(k+1)}\,.
	\end{equation}
	For the second of the terms on the right-hand side of \eqref{Eq:L2L2Err1}, we get by the exactness \eqref{Eq:ExactQuad} of the quadrature in time for all $p\in \mathbb P_{2k}(I_n;\R)$ along with \eqref{Eq:ErrEstSTnc01} that 
	\begin{equation}
		\label{Eq:L2L2Err3}
		\| I_\tau U - U_{\tau,h}\|_{\nu}^2 = \| I_\tau U - U_{\tau,h}\|_{\tau,\nu}^2 = \| U - U_{\tau,h}\|_{\tau,\nu}^2 \leq  C  (1+T) \operatorname{e}^{2\nu T} (\tau^{2(k+1)} + T h^{2r})\,.
	\end{equation}
	Together, \eqref{Eq:L2L2Err1} to \eqref{Eq:L2L2Err3} prove the assertation \eqref{Eq:L2L2Err0} . 
\end{mproof}

By Cor.\ \ref{Lem:L2L2Err}, the key error estimate \eqref{KeyRes} of this work is thus proved.

\section{Summary and outlook}
\label{Sec:SumOut}

In this work we presented and analyzed the numerical approximation of a prototype hyperbolic-parabolic model of dynamic poro- or thermoelasticity that is rewritten as a first-order evolutionary system in space and time such that the solution theory of Picard \cite{P09,STW22} becomes applicable. A family of discontinuous Galerkin (DG) schemes in space and time was studied where the innovation came through the discontinuous Galerkin discretization in space of the first-order formulation. By a consistent definition of the first-order spatial differential operators on broken polynomials spaces and the addition of boundary correction terms the mathematical evolutionary structure of the continuous problem was preserved on the fully discrete level. Well-posedness of the fully discrete system and error estimates were proved. The numerial evaluation of the approach and computational studies with comparison to the three-field formulation of \cite{ABMS23_2,ABMS23,BKR22} remain a work for the future. Further, error control in higher order norms (like the usual DG norm, cf.\ \cite{PE12}) involving the broken gradient is of interest and remains a future task. The optimality of the error estimates \eqref{Eq:ErrEstSTnc01} and \eqref{Eq:L2L2Err0} with respect to the rate of convergence in space still needs further elucidation. An improvement to convergence of order $r+1$ in space might become feasible, which is shown in \cite{B23} for an equal-order approximation of the second-order in space displacement presssure formulation of \eqref{Eq:HPS}. However, the coupling mechanisms of the unknowns in the model equations and the abstract evolutionary Problem~\ref{Prob:EP} for the holistic vector $U$ of variables do not allow such an improvement in a straightforward and obvious manner. In our error analysis, the appearance of the interpolation error $\| A_h Z\|_{\tau,\nu}$, involving first order derivatives, leads to an order reduction in space.  For this, we refer also to the results in \cite{BKR22} for the three-field formulation that also lack from optimality of the theoretical convergence rate in space, even though the latter is observed in numerical experiments.

\renewcommand{\thesection}{\Alph{section}}
\setcounter{section}{0}

\section*{Appendix}

\section{Proof of Lem.\ \ref{Lem:ErrDisc}}	

To keep this work self-contained, we present the proof of Lem.~\ref{Lem:ErrDisc}.

\begin{mproof} 
Let $U_{\tau,h}\in Y_\tau^k(H_h^{\text{dg}})$ be the solution of Problem \ref{Prob:FdPnc} and $I_\tau \Pi_h U\in  Y_\tau^k(H_h^{\text{dg}})$ its approximation in $Y_\tau^k(H_h^{\text{dg}})$ by combined interpolation and projection.  Under the regularity assumption \eqref{Eq:Reg0} there holds for the solution $U\in H_\nu(\R;D(A))$ of Problem~\ref{Prob:EP} that 
\begin{equation*}
		(\partial_t M_0 + M_1 + A)U(t) = F(t)\,,	\quad \text{for}\; t\in [0,T]\,,
\end{equation*}
such that, for $n\in \{1,\ldots,N\}$,
\begin{equation}
\label{Eq:ErrEq01}
\begin{aligned}
Q_n[(\partial_t M_0 + M_1 + A_h)U,V_{\tau,h}]_\nu & + Q_n[J_\partial(U,V_{\tau,h})+ J_\gamma(U,V_{\tau,h})]_\nu\\[1ex]
&  + \langle M_0\lsem U \rsem_{n-1} , V_{\tau,h}^{+,n-1}\rangle = Q_n[I_\tau  \Pi_h F,V_{\tau,h}]
\end{aligned}
\end{equation}
for  all $V_{\tau,h}\in Y_\tau^k(H_h^{\text{dg}})$. In \eqref{Eq:ErrEq01}, we let $A_h$ be the natural extension of Def.~\ref{Def:Ah} to $D(A)$ with
\begin{align*}
\langle \dggrd y, v_h \rangle & := \langle \grd y, v_h \rangle\,, && y\in H^1_0(\Omega)\,, \quad v_h\in (Y_h^r)^d\,,\\[1ex]
 \langle \dgdiv z, y_h \rangle &  := \langle \dv z, y_h \rangle\ - \sum_{e\in \mathcal E_h^\partial} \langle z\cdot n_e,y_h \rangle_e\,, &&z\in D(\operatorname{div})\,, \quad v_h\in Y_h^r\,,\\[1ex]
\langle \dgGrd y, v_h \rangle & := \langle \Grd y, v_h \rangle\,, && y\in H^1_0(\Omega)^d\,, \quad v_h\in (Y_h^r)^{d\times d}\cap L^2(\Omega)_{\operatorname{sym}}^{d\times d}\,,\\[1ex]
\langle \dgDiv z, y_h \rangle &  := \langle \Dv z, y_h \rangle\ - \sum_{e\in \mathcal E_h^\partial} \langle z\cdot n_e,y_h \rangle_e\,, &&z\in D(\operatorname{div})\cap L^2(\Omega)_{\operatorname{sym}}^{d\times d}\,, \quad v_h\in (Y_h^r)^d\,.
\end{align*}	
For \eqref{Eq:ErrEq01}, we note that by \eqref{Eq:DefAh2} and \eqref{Eq:FdPnc4} there holds for $V_h\in H_h^{\operatorname{dg}}$ that
\begin{align*}
& \langle A_h U, V_h \rangle + J_\partial (U, V_h ) =  - \langle \dgDiv U^{2}, V_h^{(1)}\rangle  - \langle \dgGrd U^{(1)}, V_h^{(2)}\rangle\\[0.5ex]
& \quad  +  \langle \dgdiv U^{(4)}, V_h^{(3)}\rangle  + \langle \dggrd U^{(3)}, V_h^{(4)}\rangle - \sum_{e\in \mathcal E_h^\partial} \langle U^{(2)} \cdot n_e, V_{\tau,h}^{(1)}  \rangle_e +  \sum_{e\in \mathcal E_h^\partial} \langle U^{(4)}\cdot n_e , V_{\tau,h}^{(3)}\rangle_e\\[1ex]
& = - \langle \Dv U^{2}, V_h^{(1)}\rangle  - \langle \Grd U^{(1)}, V_h^{(2)}\rangle +  \langle \dv U^{(4)}, V_h^{(3)}\rangle  + \langle \grd U^{(3)}, V_h^{(4)}\rangle = \langle A U,V_h\rangle. 
\end{align*}
Further,  we have that $J_\gamma(U,V_{\tau,h})=0$ for $U\in D(A)$. Under the assumption \eqref{Eq:Reg0}, $\lsem U  \rsem_{n-1} = 0$ is satisfied. The identity $Q_n[F,V_{\tau,h}] =Q_n[I_\tau  \Pi_h F,V_{\tau,h}]$ directly follows from \eqref{Eq:GRF}, \eqref{Def:LIO} and \eqref{Def:Pih}. Substracting now \eqref{Eq:FdPnc1} from \eqref{Eq:ErrEq01} yields with the splitting \eqref{Eq:ErrEstSTnc02} the error equation
	\begin{equation}
		\label{Eq:ErrEq02}
		\begin{aligned}
			& Q_n[(\partial_t M_0 + M_1  + A_h)E_{\tau,h},V_{\tau,h}]_\nu  + Q_n[J_\partial(E_{\tau,h},V_{\tau,h})+ J_\gamma(E_{\tau,h},V_{\tau,h})]_\nu + \langle M_0 \lsem E_{\tau,h}\rsem_{n-1}, V_{\tau,h}^{+,n-1}\rangle \\[1ex]
			& = - Q_n[(\partial_t M_0 + M_1 + A_h) Z ,V_{\tau,h}]_\nu - Q_n[J_\partial(Z,V_{\tau,h})+ J_\gamma(Z,V_{\tau,h})]_\nu- \langle M_0 \lsem Z  \rsem_{n-1}, V_{\tau,h}^{+,n-1}\rangle
		\end{aligned}
	\end{equation}
	for  all $V_{\tau,h}\in Y_\tau^k(H_h^{\text{dg}})$ and $n=1,\ldots,N$. Choosing $V_{\tau,h} = E_{\tau,h}$ and recalling that 
	\begin{equation*}
		\langle A_h E_{\tau,h},E_{\tau,h}\rangle + J_\partial(E_{\tau,h},E_{\tau,h}) = 0\,,
	\end{equation*}
	by the arguments of \eqref{Eq:Uniq05}, we get that 
	\begin{equation}
		\label{Eq:ErrEq03}
		\begin{aligned}
			& Q_n[(\partial_t M_0  + M_1 )E_{\tau,h},E_{\tau,h}]_\nu  + Q_n[J_\gamma(E_{\tau,h},E_{\tau,h})]_\nu+ \langle M_0 \lsem E_{\tau,h}\rsem_{n-1}, E_{\tau,h}^{+,n-1}\rangle \\[1ex]
			& = - Q_n[(\partial_t M_0 + M_1 + A_h) Z ,E_{\tau,h}]_\nu - Q_n[J_\partial(Z,E_{\tau,h})+ J_\gamma(Z,E_{\tau,h})]_\nu - \langle M_0 \lsem Z \rsem_{n-1}, E_{\tau,h}^{+,n-1}\rangle \\[1ex]
			& =: E_i^n
		\end{aligned}
	\end{equation}
	for $n=1,\ldots N$ and  $E_{\tau,h}^{+,n-1}:= E_{\tau,h}^+(t_{n-1})$; cf.~\eqref{Eq:DefLim}.
	
By \cite[Lem.~3.5]{FTW19}, for the left-hand side of \eqref{Eq:ErrEq03} there holds that
\begin{equation}
		\label{Eq:ErrEq04}
		\begin{aligned}
			& Q_n[(\partial_t M_0  + M_1 )E_{\tau,h},E_{\tau,h}]_\nu  + Q_n[J_\gamma(E_{\tau,h},E_{\tau,h})]_\nu + \langle M_0 \lsem E_{\tau,h}\rsem_{n-1}, E_{\tau,h}^{+,n-1}\rangle\\[1ex]  
			&\quad  \geq \gamma \| E_{\tau,h}\|_{\tau,\nu,n}^2  + |J_\gamma(E_{\tau,h},E_{\tau,h})|_{\tau,\nu,n}^2 + \frac{1}{2} \Big[ \langle M_0 E_{\tau,h}^-(t_n), E_{\tau,h}^-(t_n)\rangle  \operatorname{e}^{-2\nu \tau_n}\\[1ex]
			& \qquad -  \langle M_0 E_{\tau,h}^-(t_{n-1}), E_{\tau,h}^-(t_{n-1})\rangle + \langle M_0 \lsem E_{\tau,h} \rsem_{n-1} , \lsem E_{\tau,h}\rsem_{n-1} \rangle \Big]
		\end{aligned}
\end{equation}
for $n=1,\ldots,N$, where $E_{\tau,h}^-(t_0)= \Pi_h U_0 - U_{0,h}$. Multiplying \eqref{Eq:ErrEq04} with the weight  $\operatorname{e}^{-2\nu t_{n-1}}$, combining this with \eqref{Eq:ErrEq03}, summing up the resulting equation and neglecting the positive jump terms yield that 
	\begin{equation}
		\label{Eq:ErrEq05}
		\begin{aligned}
			\langle M_0 E_{\tau,h}^-(t_{N}), E_{\tau,h}^-(t_{N})\rangle & \operatorname{e}^{-2\nu T}  + \gamma  \| E_{\tau,h}\|_{\tau,\nu}^2  + |J_\gamma(E_{\tau,h},E_{\tau,h})|_{\tau,\nu}^2  \\[1ex]
			& \leq C \Big(\langle M_0 E_{\tau,h}^-(t_0), E_{\tau,h}^-(t_0)\rangle + \sum_{n=1}^N \operatorname{e}^{-2\nu t_{n-1}} E_i^n \Big)\,.
		\end{aligned}
	\end{equation}

To bound $E_i^n$ in \eqref{Eq:ErrEq05}, we need auxiliary results. Recalling the exactness of the quadrature formula \eqref{Eq:GRF} for all $w\in P_{2k}(I_n;\R)$ along with \eqref{Def:IPPk} and using integration by part, we get that
	\begin{equation*}
		\label{Eq:ErrEq06}
		\begin{aligned}
			& Q_n[\partial_t M_0 I_\tau \Pi_h U,E_{\tau,h}]_\nu   = \langle \partial_t M_0 I_\tau \Pi_h U ,E_{\tau,h}\rangle_{\nu,n} \\[1ex]
			& \qquad = \underbrace{\langle \operatorname{e}^{-2\nu(t-t_{n-1})}M_0 I_\tau\Pi_h  U,E_{\tau,h}\rangle_H\Big|_{t_{n-1}}^{t_n}}_{=:a} - \langle M_0 I_\tau \Pi_h  U, \partial_t E_{\tau,h}\rangle_{\nu,n}\\[1ex]
			& \qquad  = a - Q_n [M_0 I_\tau \Pi_h U,\partial_t E_{\tau,h}]_\nu = a - Q_n [M_0 \hItau \Pi_h U,\partial_t E_{\tau,h}]_\nu\\[1ex] 
			& \qquad = a - \langle M_0 \hItau \Pi_h U,\partial_t E_{\tau,h} \rangle_{\nu,n}\\[1ex]
			& \qquad = a + \langle \partial_t M_0 \hItau\Pi_h  U, E_{\tau,h} \rangle_{\nu,n} - \underbrace{\langle \operatorname{e}^{-2\nu(t-t_{n-1})} M_0 \hItau \Pi_h U, E_{\tau,h} \rangle_H\Big|_{t_{n-1}}^{t_n}}_{=:b}\\[1ex]
			& \qquad  = a-b + Q_n[\partial_t M_0 \hItau \Pi_h U, E_{\tau,h}]_\nu\,.
		\end{aligned}
	\end{equation*}
From this along with the definition of $I_\tau$ and $\hItau$ in \eqref{Def:LIO} and \eqref{Def:LIOc}, respectively, we conclude that 
	\begin{equation}
		\label{Eq:ErrEq07}
		\begin{aligned}
			Q_n[\partial_t M_0 I_\tau \Pi_h  U,E_{\tau,h}]_\nu + \langle M_0 \lsem I_\tau \Pi_h  U \rsem_{n-1}, E_{\tau,h}^{+,n-1}\rangle_H = Q_n[\partial_t M_0 \hItau \Pi_h  U, E_{\tau,h}]_\nu\,.
		\end{aligned}
	\end{equation}
By the inequalities of Cauchy--Schwarz and Cauchy--Young  there holds that (cf.\ \eqref{Eq:FdPnc3} and \eqref{Eq:ErrDisc00})
\begin{equation}
	\label{Eq:ErrEq09}
	|J_\partial(Z, E_{\tau,h})| \leq C  J_\partial^n (Z,Z) + \beta J_\gamma (E_{\tau,h},E_{\tau,h})
\end{equation}
for any $\beta >0$. Then, for the interpolation error $E_i^n$ defined in \eqref{Eq:ErrEq05}  we obtain by \eqref{Eq:ErrEq07}, \eqref{Eq:ErrEq09}  and the inequalities of Cauchy--Schwarz and Cauchy--Young that 
	\begin{equation}
		\label{Eq:ErrEq8}
		\begin{aligned}
			\sum_{n=1}^N & \operatorname{e}^{-2\nu t_{n-1}}  | E_i^n|  \leq C \Big( \| \partial_t M_0 (U - \hItau  \Pi_h U) \|_{\tau,\nu}^2 +  \| M_1  Z \|_{\tau,\nu}^2 + \| A_h Z \|_{\tau,\nu}^2 +  | J_\partial^n (Z,Z)|_{\tau,\nu}^2 \\[1ex]
			&  \quad +  | J_\gamma (Z,Z)|_{\tau,\nu}^2 +T\max_{1\leq n\leq N} \Big\{ \|M_0 (\Pi_h U^+(t_{n-1}) - I_\tau \Pi_h U^+(t_{n-1}))\|_H^2 \operatorname{e}^{-2\nu t_{n-1}} \Big\} \Big)\\[1ex]
			& \quad + \beta_1 \| E_{\tau,h}\|_{\tau,\nu}^2 + \beta_2 | J_\gamma (E_{\tau,h},E_{\tau,h})|_{\tau,\nu}^2
		\end{aligned}
	\end{equation}
	for any $\beta_1,\beta_2 >0$. Finally, combining \eqref{Eq:ErrEq05} with \eqref{Eq:ErrEq8} and choosing $\beta_1$ and $\beta_2$ sufficiently small, proves the assertion \eqref{Eq:ErrDisc0}.
\end{mproof}

\section{Proof of Lem.\ \ref{Lem:ErrDiscInf}}

To keep this work self-contained, we present the proof of Lem.~\ref{Lem:ErrDiscInf}.

\begin{mproof}
	Using an idea of \cite[Cor.\ 2.1]{AM04}, for $E_{\tau,h}\in Y_\tau^k(H_h^{\text{dg}})$ we define the local interpolant
	\begin{equation*}
		\widehat E_{\tau,h}:= \InGR \Phi\,, \quad \text{with } \;\; \Phi := \frac{\tau_n}{t - t_{n-1}} E_{\tau,h}\,, \quad \text{for }\; t\in I_n \,, \; n=1,\ldots,N\,,
	\end{equation*}
	where the local Lagrange interpolation operator $\InGR: C(I_n;B) \to \mathbb P_k (I_n;B)$, for $n\in \{1,\ldots,N\}$,  satisfies 
	\begin{equation*}
		\InGR f(t_{n,\mu}) = f (t_{n,\mu})\,, \quad \text{for }\; \mu = 0,\ldots ,k\,,
	\end{equation*}
	for the quadrature nodes $t_{n,\mu}\in I_n$, for $\mu = 0,\ldots, k$, of the (non-weighted) Gauss--Radau formula on $I_n$. Then, there holds that  	
	\begin{equation}
		\label{Eq:ErrDiscInf1}
		\begin{aligned}
			\langle M_0 \widehat E_{\tau,h}(t_{n,\mu}),\widehat E_{\tau,h}(t_{n,\mu})\rangle_H & = \frac{\tau_n^2}{(t_{n,\mu}-t_{n-1})^2}\langle M_0 E_{\tau,h}(t_{n,\mu}), E_{\tau,h}(t_{n,\mu})\rangle_H \\[1ex]  
			& \geq  \langle M_0 E_{\tau,h}(t_{n,\mu}), E_{\tau,h}(t_{n,\mu})\rangle_H \,.
		\end{aligned}
	\end{equation}
	By \cite[Lem.~3.10]{FTW19}, that is based on \cite[Lem.~2.1]{AM04}, along with \eqref{Eq:ErrDiscInf1} we obtain that 
	\begin{equation}
		\label{Eq:ErrDiscInf2}
		\begin{aligned}
			Q_n[\partial_t M_0 E_{\tau,h}, 2 \widehat E_{\tau,h}]_\nu + \langle M_0 E_{\tau,h}^+(t_{n-1}), 2\widehat E_{\tau,h}^+(t_{n-1}) \rangle_H & \geq \dfrac{1}{\tau_n} Q_n [M_0 \widehat E_{\tau,h},\widehat E_{\tau,h}]_\nu\\[1ex] 
			& \geq \dfrac{1}{\tau_n} Q_n [M_0 E_{\tau,h}, E_{\tau,h}]_\nu\,.
		\end{aligned}
	\end{equation}
	By the norm equivalence 
	\begin{equation*}
		\sup_{t\in [0,t]} |w(t))| \leq C_e \|w \|_{L^1((0,1);\R)}\,, \quad \text{for }\; w\in \mathbb P_k ([0,1];\R)\,, 
	\end{equation*}
	along with the transformation of $[t_{n-1},t_n]$ to $[0,1]$ we have that 
	\begin{equation}
		\label{Eq:ErrDiscInf3}
		\sup_{t\in I_m}  \langle M_0 E_{\tau,h}(t), E_{\tau,h}(t)\rangle_H  \leq \frac{C_e}{\tau_n} \operatorname{e}^{2\nu \tau_n} Q_n[M_0 E_{\tau,h},E_{\tau,h}]_\nu \leq \frac{C}{\tau_n} Q_n[M_0 E_{\tau,h} E_{\tau,h}]_\nu
	\end{equation}
	with $C:= C_e \operatorname{e}^{2\nu T}\geq \max_{n=1,\ldots,N}\{\operatorname{e}^{2\nu \tau_n}\} C_e$. Further, by \eqref{Eq:DefAh3} along with \eqref{Eq:FdPnc4} we have that 
	\begin{equation}
		\label{Eq:ErrDiscInf4}
		\begin{aligned}
		& Q_n[A_h E_{\tau,h},2 \widehat E_{\tau,h}]_\nu + Q_n[J_\partial (E_{\tau,h},2 \widehat E_{\tau,h})]_\nu  \\[1ex] 
		& = \dfrac{\tau_n}{2} \sum_{\mu = 0}^k \hat \omega_\mu \dfrac{2\tau_n}{t_{n,\mu}-t_{n-1}} \Big(\langle A_h E_{\tau,h}(t_{n,\mu}),E_{\tau,h}(t_{n,\mu}) \rangle_H + J_\partial(E_{\tau,h}(t_{n,\mu}),E_{\tau,h}(t_{n,\mu}) \Big)\\[1ex]
		& = 0\,.
		\end{aligned}
	\end{equation}
	Combining \eqref{Eq:ErrDiscInf3} with \eqref{Eq:ErrDiscInf2} and then using \eqref{Eq:ErrDiscInf4}, it follows that 
	\begin{equation}
		\label{Eq:ErrDiscInf5}
		\begin{aligned}
			 \sup_{t\in I_m} & \langle M_0 E_{\tau,h}(t),  E_{\tau,h}(t)\rangle_H\\[1ex] 
			 & \leq C \Big(Q_n[\partial_t M_0 +M_1 + A_h) E_{\tau,h}, 2 \widehat E_{\tau,h}]_\nu + Q_n[J_\partial (E_{\tau,h},2 \widehat E_{\tau,h})+J_\gamma (E_{\tau,h},2 \widehat E_{\tau,h})]_\nu \\[1ex]
			& \quad  + \langle M_0 \lsem E_{\tau,h}\rsem_{n-1} , 2\widehat E_{\tau,h}^+(t_{n-1}) \rangle_H  - Q_n[M_1 E_{\tau,h}, 2 \widehat E_{\tau,h}]_\nu  \\[1ex]
			& \quad - Q_n[J_\gamma (E_{\tau,h},2 \widehat E_{\tau,h})]_\nu + \langle M_0 E_{\tau,h}^-(t_{n-1}), 2\widehat E_{\tau,h}^+(t_{n-1}) \rangle_H\Big) \,.
		\end{aligned}
	\end{equation}
	Using the error equation \eqref{Eq:ErrEq02} with test function $V_{\tau,h}=2\widehat E_{\tau,h}$, we deduce from \eqref{Eq:ErrDiscInf5} that 
	\begin{equation}
		\label{Eq:ErrDiscInf6}
		\begin{aligned}
			\sup_{t\in I_m}  & \langle M_0 E_{\tau,h}(t),   E_{\tau,h}(t)\rangle_H\\[1ex] 
			&  \leq C \Big(-Q_n[\partial_t M_0 +M_1 + A_h) Z, 2 \widehat E_{\tau,h}]_\nu - Q_n[J_\partial (Z,2 \widehat E_{\tau,h})+J_\gamma (Z,2 \widehat E_{\tau,h})]_\nu \\[1ex]
			& \quad  - \langle M_0 \lsem Z \rsem_{n-1} , 2\widehat E_{\tau,h}^+(t_{n-1}) \rangle_H - Q_n[M_1 E_{\tau,h}, 2 \widehat E_{\tau,h}]_\nu \\[1ex] 
			& \quad - Q_n[J_\gamma (E_{\tau,h},2 \widehat E_{\tau,h})]_\nu + \langle M_0 E_{\tau,h}^-(t_{n-1}), 2\widehat E_{\tau,h}^+(t_{n-1}) \rangle_H\Big) \,.
		\end{aligned}
	\end{equation}
	
	Next, we bound the right-hand side in \eqref{Eq:ErrDiscInf6}. For this, we use the boundedness of $M_1$, the identity $\langle M_0 Z^-(t_{n-1}),\widehat E_{\tau,h}^+(t_{n-1}) \rangle_H = 0$ by the definition of $Z$ in \eqref{Eq:ErrEstSTnc02} and $I_\tau$ in \eqref{Def:LIO} as well as
	\begin{equation*}
		\langle M_0 u,v\rangle_H = \langle M_0^{1/2} u,M_0^{1/2}v\rangle_H \leq \langle M_0^{1/2} u,M_0^{1/2}u\rangle_H \langle M_0^{1/2} v,M_0^{1/2}v\rangle_H \,, \quad \text{for }\; u,v\in H\,, 
	\end{equation*}
	by the non-negativity and selfadjointness of $M_0$. Similarly to \eqref{Eq:ErrEq8}, we then get from \eqref{Eq:ErrDiscInf6} that 
	\begin{equation}
		\label{Eq:ErrDiscInf7}
		\begin{aligned}
			\sup_{t\in I_n} & \langle M_0 E_{\tau,h}(t),  E_{\tau,h}(t)\rangle_H\\[1ex]
			&   \leq C \Big( \| \partial_t M_0 (U - \hItau  \Pi_h U) \|_{\tau,\nu,n}^2 +  \| M_1  Z \|_{\tau,\nu,n}^2 + \| A_h Z \|_{\tau,\nu,n}^2   \\[1ex]
		   & +  |J^n_\partial (Z,Z)|_{\tau,\nu,n} + |J_\partial (Z,Z)|_{\tau,\nu,n} + \|M_0 Z^+(t_{n-1})\|_H^2\Big)\\[1ex]	   
		   &  + \underbrace{\alpha_1\langle M_0 E_{\tau,h}^-(t_{n-1}), E_{\tau,h}^-(t_{n-1})\rangle_H + \alpha_2 \|M_1 \|^2 \|E_{\tau,h}\|_{\tau,\nu,n}^2+ \alpha_3 |J_\gamma(E_{\tau,h},E_{\tau,h})|_{\tau,\nu,n}}_{=:G_n(E_{\tau,h})} \\[1ex]
			& + \beta_1  \langle 2 M_0 \widehat E_{\tau,h}^+(t_{n-1}), 2\widehat E_{\tau,h}^+(t_{n-1}) \rangle_H+ \beta_2 Q_n[2\widehat E_{\tau,h},2\widehat E_{\tau,h}]_\nu + \beta_3 Q_n[J_\gamma(2\widehat E_{\tau,h},2\widehat E_{\tau,h})]_\nu
		\end{aligned}
	\end{equation}
	with some constants $\alpha_i,\beta_i> 0$, for $i=1,2,3$. The term $G_n(E_{\tau,h})$ is bounded by means of \eqref{Eq:ErrDisc0}, that still holds if $E_{\tau,h}^-(t_{N})$ is replaced by $E_{\tau,h}^-(t_{n})$ for $n=1,\ldots,N-1$. Further, noting that (cf.~\cite[Cor.~1.5]{TW16}) 
	\begin{equation*}
		\frac{\tau_n}{t_{n,\mu}-t_{n-1}} \leq \frac{\tau_n}{t_{n,0}-t_{n-1}} \leq \frac{1}{\delta} \,, \quad \text{for }\; n\in \{1,\ldots, N\}\,,
	\end{equation*}
	for some $\delta >0$ depending on $\nu$ and $T$ only, we have that 
	\begin{subequations}
		\label{Eq:ErrDiscInf8}
		\begin{alignat}{2}
			Q_n[\widehat E_{\tau,h},\widehat E_{\tau,h}]_\nu &\leq \dfrac{1}{\delta^2} \|E_{\tau,h} \|^2_{\tau,\nu,n}\,,\\[1ex]
			Q_n[J_\gamma(\widehat E_{\tau,h},\widehat E_{\tau,h})]_\nu & \leq \dfrac{1}{\delta^2} |J_\gamma(E_{\tau,h},E_{\tau,h})|_{\tau,\nu,n}^2  \\[1ex]
			\langle M_0 \widehat E_{\tau,h}^+(t_{n-1}), \widehat E_{\tau,h}^+(t_{n-1}) \rangle_H & \leq \dfrac{1}{\delta^2}  \sup_{t\in I_n} \, \langle M_0 E_{\tau,h}(t), E_{\tau,h}(t)\rangle_H\,.
		\end{alignat}
	\end{subequations}	
	Finally, combining \eqref{Eq:ErrDiscInf7} for a sufficiently small choice of $\beta_1$ to $\beta_3$ with \eqref{Eq:ErrDiscInf8} and using \eqref{Eq:ErrDisc0} proves the assertion \eqref{Eq:ErrDiscInf0}. 
\end{mproof}

\end{document}